


\def\2{{1\over 2}}

\def\d{\delta}
\def\a{\alpha}
\def\b{\beta}
\def\g{\gamma}

\def\s{\sigma}
\def\e{\epsilon}
\def\l{\lambda}

\def\fun#1#2#3{#1\colon #2\rightarrow #3}

\def\frac#1#2{{{#1} \over {#2}}}

\def\sqr{\sqrt}
\def\st{\;\colon\;}
\def\tends{\rightarrow}

\def\dr{ {\rm d} }

\def\R{{\bf R}}
\def\N{{\bf N}}

\def\thm#1{\vskip 1 pc\noindent{\bf Theorem #1.\quad}\sl}
\def\lem#1{\vskip 1 pc\noindent{\bf Lemma #1.\quad}\sl}

\def\proof{\rm\vskip 1 pc\noindent{\bf Proof.\quad}}
\def\fin{\par\hfill $\backslash\backslash\backslash$\vskip 1 pc}
\def\txt#1{\quad\hbox{#1}\quad}

\def\L{{\cal L}}

\def\s{\sigma}

\def\2{\frac{1}{2}}
\def\inn#1#2{{\langle #1 ,#2\rangle}}

\def\diff{\frac{\dr}{\dr t}}

\def\part{{\partial_{x}}}

\def\tr{{{}^{t}}}

\def\ec{{\cal E}}

\def\dc{{\cal D}}
\def\oc{{\cal O}}
\def\pc{{\cal P}}

\def\mc{{\cal M}}

\def\rc{{\cal R}}

\def\listo{{i_1\dots i_{l}}}



\baselineskip= 17.2pt plus 0.6pt
\font\titlefont=cmr17
\centerline{\titlefont Harmonic embeddings }
\vskip 1 pc
\centerline{\titlefont of the Stretched Sierpinski Gasket}
\vskip 4pc
\font\titlefont=cmr12
\centerline{         \titlefont {Ugo Bessi}\footnote*{{\rm 
Dipartimento di Matematica, Universit\`a\ Roma Tre, Largo S. 
Leonardo Murialdo, 00146 Roma, Italy.}}   }{}\footnote{}{
{{\tt email:} {\tt bessi@matrm3.mat.uniroma3.it} Work partially supported by the PRIN2009 grant "Critical Point Theory and Perturbative Methods for Nonlinear Differential Equations}} 
\vskip 0.5 pc
 
\par
\vskip 2pc
\centerline{\bf Abstract} 

P. Alonso-Ruiz, U. Freiberg and J. Kigami have defined a large family of resistance forms on the Stretched Sierpinski Gasket $G$. In the present paper we introduce a system of coordinates on $G$ (technically, an embedding of $G$ into $\R^2$) such that 

\noindent$\bullet$) these forms are defined on $C^1(\R^2,\R)$ and 

\noindent$\bullet$) all affine functions are harmonic for them. 

We do this adapting a standard method from the Harmonic Sierpinski Gasket: we start finding a sequence $G_l$ of pre-fractals such that all affine functions are harmonic on 
$G_l$. After showing that this property is inherited by the stretched harmonic gasket $G$, we use the formula for the Laplacian of a composition to prove that, for a natural measure 
$\mu$ on $G$, $C^2(\R^2,\R)\subset\dc(\Delta)$ and Teplyaev's formula for the Laplacian of $C^2$ functions holds. Lastly, we use the expression for $\Delta u$ to show that the form we have found is closable in 
$L^2(G,\mu)$.

\vskip 2 pc
\centerline{\bf  Introduction}
\vskip 1 pc

In [1], P.Alonso-Ruiz, U. Freiberg and J. Kigami extend a result of [2] and define a family of resistance forms on the Stretched Sierpinski Gasket; we briefly sketch the construction of this relative of the Sierpinski Gasket. 

One considers the vertices of an equilateral triangle in $\R^2$; 
$$A=\left(
\matrix{
0\cr
0
}
\right)  ,  \qquad
B=\left(
\matrix{
\frac{\sqrt 3}{2}\cr
\2
}
\right)  , \qquad
C=\left(
\matrix{
\frac{\sqrt 3}{2}\cr
-\2
}
\right)   \eqno (1)$$
and the three contractions, depending on a parameter $\a\in(0,\2)$:
$$\psi_1\left(
\matrix{
x\cr
y
}
\right)=
\a\left(
\matrix{
x\cr
y
}
\right),\qquad
\psi_2\left(
\matrix{
x\cr
y
}
\right)=\a\left[
\left(
\matrix{
x\cr
y
}
\right)-B
\right] +B,$$
$$\psi_2\left(
\matrix{
x\cr
y
}
\right)=\a\left[
\left(
\matrix{
x\cr
y
}
\right)-C
\right] +C  .  $$
Denoting by $PQ$ the segment joining $P,Q\in\R^2$, the Stretched Sierpinski Gasket is the unique compact set $G\subset\R^2$ such that
$$G=\bigcup_{i=1}^3
\psi_i(G)\cup
[
\psi_1(B)\psi_2(A)\cup\psi_1(C)\psi_3(A)\cup\psi_2(C)\psi_3(B)
]  .  \eqno (2)$$
The segments in figure 1 below are the first stage of the construction of $G$ or, technically speaking, the first pre-fractal $G_1$ of formula (4) below; in the figure the maps $\psi_i$ are called $F_i$ because this is the notation we'll adopt when the maps are affine, but not necessarily homotheties.

\centerline{Figure 1}

\vskip 1pc

Many fractals $G$ are induced by Iterated Function Systems: roughly speaking, there are invertible contractions $\fun{F_1,\dots,F_p}{\R^2}{\R^2}$ such that a "homogeneous" version of (2) holds: 
$$G=\bigcup_{i=1}^p F_i(G)  .  $$
For this kind of sets, the construction of a resistance form is a staple of fractal theory: we refer the reader to [10] and [15] for a purely analytical construction; [4] has a different approach, which uses the dynamics on $G$ (as we shall see in section 1 below, this dynamics exists if the maps $F_i$ are invertible) and the fact that the "carr\'e\ de champs" is closely related to a Gibbs measure. 

This method cannot be applied directly to the stretched gasket because the maps of [2] and [3] are not invertible: to build a form on the stretched gasket one has to stretch the approach of [10]. In the present paper we re-prove the result of [1] using the dynamical approach of [4]; a side result will be a harmonic embedding of $G$ in $\R^2$, by which we mean that all affine functions are harmonic on the image of $G$. The name comes from the analogous property of the harmonic Sierpinski gasket (see [9] for a thorough study of this set). 

We briefly explain our construction. A naive idea would be to replace the homotheties 
$\psi_i$ with affine maps $F_i$ such that all pre-fractals are harmonic: by this we mean that all affine functions are harmonic on them. This does not work for two reasons: 

\noindent 1) it doesn't yield a natural way to insert infinitely many parameters, as in [1].

\noindent 2) The "cable" part of the Dirichlet form does not converge (a fact which will be evident in section 5). The term "cable" comes from [1] and we shall define it precisely in formula (6) below. 

Thus, we follow a different path: we consider a sequence 
$\{ (F^i_1,F^i_2,F^i_3) \}_{i\ge 1}$ of triples of affine maps. The maps of each triple are tied by a very strong symmetry property and the first elements of the triples have the following form, strongly reminiscent of the harmonic gasket: for $\e_i\in(0,1)$, 
$$F^i_1\left(
\matrix{
x\cr
y
}
\right)  =
\e_i\left(
\matrix{
\frac{3}{5}& 0\cr
0&\frac{1}{5}
}
\right)  .  \eqno (3)$$

Now we define the $l$-th pre-fractal. We set 
$F_{i_1\dots i_s}=F_{i_1}^1\circ\dots\circ F_{i_s}^s$ with the convention that 
$F_{i_1\dots i_s}$ is the identity when $s=0$; we define 
$$G_0=\tilde G_0=AB\cup BC\cup AC   .   $$
For $l\ge 0$ we set 
$$G_l=\tilde G_l\cup C_l    \eqno (4)$$
where 
$$\tilde G_l\colon=\bigcup_{i_1\dots i_l\in\{ 1,2,3 \}}
F_{i_1\dots i_l}(G_0)  .  \eqno (5)$$
As we shall see, $\tilde G_l$ is a union of disjoint triangles; the number $\e_l$ in (3) is the {\it relative} size of $F_{i_1\dots i_l}(G_0)$ with respect to $F_{i_1\dots i_{l-1}}(G_0)$ into which it is contained; this fact somehow recalls the construction of the Cantor set of positive measure ([18]). The set $C_l$ of (4) is the set of the cables connecting the triangles: namely, we set 
$$\Omega^s=\{
F_1^s(B)F_2^s(A),F_1^s(C)F_3^s(A),F_2^s(C)F_3^s(B)
\}     \eqno (6)$$
and 
$$C_l=\bigcup_{s=1}^l\bigcup_{i_1\dots i_{s-1}\in\{ 1,2,3 \}}
F_{i_1\dots i_{s-1}}(\Omega^s)  .  $$
With the convention we have adopted, when $s=1$ 
$F_{i_1\dots i_{s-1}}(\Omega^s)=\Omega^1$ is the set of the first generation cables which connect together the first generation triangles $F^1_j(G_0)$; it is the situation depicted in figure 1. When $s=2$ we get the cables connecting the second generation triangles $F_{i_1i_2}(G_0)$ and so on. The stretched gasket is the Haussdorff limit of the sets $G_l$ as $l\tends+\infty$. 

We define a Dirichlet form on $G_l$. For the sequence $\{ \e_i \}_{i\ge 1}$ of (3) we set 
$$\rc=\left\{ 
\e_i\left( \frac{3}{5},\frac{1}{5} \right) 
\right\}_{i\ge 1},\qquad \s\rc=
\left\{ \e_{i+1} \left( \frac{3}{5},\frac{1}{5} \right) 
\right\}_{i\ge 1}  .  \eqno (7)$$
If $P,Q\in\R^2$, we parametrise the segment $PQ$ in the standard way:  
$$\fun{\g_{PQ}}{[0,1]}{\R^2},\qquad
\g_{PQ}(t)=(1-t)P+tQ  .  \eqno (8)$$
For $l\ge 1$, we set 
$$\fun{\ec_{\rc,l}}{C^1(\R^2,\R)\times C^1(\R^2,\R)}{\R}$$
$$\ec_{\rc,l}(\nabla u,\nabla v)=
\ec^1_{\rc,l}(\nabla u,\nabla v)+\ec^2_{\rc,l}(\nabla u,\nabla v)  .  \eqno (9)$$
We define the two bilinear forms on the right hand side of (9) in a similar way to [1]; with the same notation as above we set 
$$\ec^1_{\rc,l}(\nabla u,\nabla v)=$$
$$\frac{a}{\left(\frac{3}{5}\right)^l\cdot\prod_{i=1}^l\e_i^2}
\sum_{\listo\in\{ 1,2,3 \}}\sum_{j\in\{ AB,BC,AC \}}
\int_0^1\diff u\circ F_\listo\circ\g_j(t)\cdot \diff v\circ F_\listo\circ\g_j(t) \dr t  \eqno (10)$$
where $a>0$ is a suitable number. When $l=1$, the integrals above are on the three small triangles of figure 1; when $l\ge 2$ we are integrating on the pre-fractal $\tilde G_l$ defined by (5). 

Integrating on the cables and their images gives us $\ec^2_l$; in other words, for a suitable $b>0$, for $\tilde\l_{\rc,l}$, $\tilde\e^l_{\rc,s}$ as in formula (5.1) below and $\Omega^s$ as in (6), we set 
$$\ec^2_{\rc,l}(\nabla u,\nabla v)=$$
$$\sum_{s=1}^{l}\sum_{i_1\dots i_s\in\{ 1,2,3 \}}\sum_{j\in\Omega^s}
\frac{b}{\tilde\l_{\rc,s-1}\cdot\tilde\e^l_{\rc,s}(1-\e_s)}\cdot
\int_0^1\diff u\circ F_{i_1\dots i_{s-1}}\circ\g_j(t)\cdot 
\diff v\circ F_{i_1\dots i_{s-1}}\circ\g_j(t)  \dr t   .  \eqno (11)$$

We shall see in section 6 below that the form $\ec_{\rc,l}$ of (9) converges to a bilinear form 
$$\fun{\ec_\rc}{C^1(\R^2,\R)\times C^1(\R^2,\R)}{\R}  .  $$

The major problem here is not convergence, which follows in a standard way from dynamical considerations: it is to show that $\ec_\rc$ has the properties of a resistance form (a definition is in [1], the theory is in [10] and [11]). We shall look at this from another perspective, that of Dirichlet forms; the relationship between resistance and Dirichlet forms is in section 8 of [11]. In this setting, the problem becomes proving that $\ec_\rc$ extends to a local, regular Dirichlet form in $L^2(G,\mu)$  for some reasonable measure 
$\mu$. The part that requires more work is closability and it is here that we need harmonicity: namely, our choice of the triples $\{ (F^i_1,F^i_2,F^i_3) \}$ implies that affine functions are harmonic for all forms $\ec_{\rc,l}$. The limit form 
$\ec_\rc$ will inherit this property, i. e. 
$$\ec_\rc(\nabla u,\nabla v)=0  \eqno (12)$$
whenever $\fun{u}{\R^2}{\R}$ is affine and the test function $v\in C^1(\R^2,\R)$ satisfies 
$v(A)=v(B)=v(C)=0$. As we shall see, (12) allows us to show that $C^2$ functions admit a Laplacian which is given by Teplyaev's formula ([19]), i. e. formula (7.1) below. The proof follows [6]: essentially, we see $u(x,y)$ as the composition of $u$ with the two harmonic functions $x$ and $y$, and then use the standard method to find the Laplacian of a composition. Once we have Teplyaev's formula, integration by parts shows easily that the map $\fun{}{u}{\ec_\rc(u,u)}$ is lower semicontinuous on $L^2(G,\mu)$; by [14], this is tantamount to the fact that $\ec$ is closable on this space. 

We briefly outline the difference between the approach of [1] and ours. In [1], the fractal is defined by the homotheties $\psi_i$ we defined after (1); the authors look for resistance forms on the pre-fractals which increase to a resistance form defined on a certain subspace of continuous functions on the fractal. In this approach, the harmonic immersion comes last, since existence of harmonic functions for resistance forms is a standard fact ([11]). 

On the other side, we are looking for triples of affine maps $(F_1^i,F_2^i,F_3^i)$ such that 

\noindent 1) the form $\ec_\rc$ has the simple domain $C^1(\R^2,\R)$. 

\noindent 2) Affine functions are harmonic on the pre-fractals and, taking limits, on the fractal. As a consequence,  

\noindent 3) for a natural measure $\mu$ on $G$, the Laplacian has a simple domain, 
$C^2(\R^2,\R)$, and a simple expression. 

As for the convergence of the forms on the pre-fractals, this will follow from standard results of Dynamical Systems which imply convergence to the Gibbs measure.

At the end, we shall have proven the following theorem. 

\thm{1} Let the numbers $\{ \e_i \}_{i\ge 1}\subset(0,1)$ in the sequence $\rc$ of (7) satisfy  
$$\prod_{i\ge 1}\e_i>0   .  $$
Then, $\rc$ determines a sequence of triples of affine maps 
$\{ (F^i_1,F^i_2,F^i_3) \}_{i\ge 1}$ such that the following points hold. 

\noindent 1) Let $\ec_{\rc,l}$ be the form on the pre-fractal $G_l$ defined by (9) and let 
$\ec_\rc$ be the form on $G$ defined by (6.10) below. Then, for all 
$u,v\in C^1(\R^2,\R)$, we have that 
$$\ec_{\rc,l}(\nabla u,\nabla v)\tends\ec_\rc(\nabla u,\nabla v)
\txt{as}l\tends+\infty  .  $$

\noindent 2) Affine functions are harmonic for $\ec_\rc$, i. e. formula (12) holds. 

\noindent 3) For the finite Borel measure $\mu$ on $G$ of formula (6.14) below, 
$C^2(\R^2,\R)$ is contained in the domain of the Laplacian. 

\noindent 4) The form $\ec_\rc$ of point 1) above extends to a regular, local Dirichlet form in $L^2(G,\mu)$. 

\noindent 5) Defining $\ec^2_{\{ \Omega^j \},s}$ as in (3.2) and $\tilde\e^\infty_{\rc,s}$ as after formula (6.3) we have that 
$$\ec_\rc(\nabla u,\nabla v)=
\frac{1}{ \e_1^2\cdot\frac{3}{5} }\sum_{i=1}^3
\ec_{\s\rc}(\nabla u\circ F^1_i,\nabla v\circ F^1_i)+
\frac{1}{\tilde\e^\infty_{\rc,1}}\ec^2_{\{ \Omega^j \},1}(\nabla u,\nabla v)  .  $$

\rm
 
\vskip 1pc

The paper is organised as follows. In section 1 below, we recall from [4] the dynamical definition of Kusuoka's measure and bilinear form. In section 2, we write down the maps $F_1^i$ and define $F_2^i$, $F_3^i$ by symmetry; in section 3 we define energy and harmonicity on the pre-fractal $G_l$; in section 4 we choose the coefficients of $F_1^i$ (which determine those of $F_2^i$ and $F_2^i$ by symmetry) so that the set $G_1$ of (4) is harmonic. In section 5 we show that the sets $G_l$ are harmonic for all $l\ge 1$. In section 6 we take limits as $l\tends+\infty$ and build the form $\ec_\rc$ and the measure 
$\mu$; in section 7, we end the proof of theorem 1, proving closability. 

\vskip 2pc

\centerline{\bf \S 1}
\centerline{\bf A dynamical approach to Kusuoka's measure and bilinear form}
\vskip 1pc

\noindent{\bf Fractal sets.} We begin considering a fractal $\tilde G$ which satisfies hypotheses (F1), (F2) and (F3) below; this is the case when the maps $\{ F^l_j \}_{j=1}^p$ do not depend on the iteration $l$; the case when they depend on the iteration is treated at the end of this section. 

\noindent {\bf(F1)} There are $p$ invertible affine maps  
$$\fun{F_1,\dots,F_p}{\R^d}{\R^d}  $$
satisfying 
$$\eta\colon=\sup_{i\in(1,\dots,p)}Lip(F_i)<1  .  \eqno (1.1)$$
By theorem 1.1.7 of [10], there is a unique non empty compact set $\tilde G\subset\R^d$ such that
$$\tilde G=\bigcup_{i=1}^p F_i(\tilde G) . \eqno (1.2)$$

If (F1) holds, then the dynamics of $F$ on $\tilde G$ can be coded. Indeed, we define 
$\Sigma$ as the space of sequences 
$$\Sigma=\{ 1,\dots,p \}^\N=
\{  \{ x_i \}_{i\ge 1}\st x_i\in(1,\dots,p),\quad\forall i\ge 1
\}    $$
with the product topology. This is a metric space; for instance, if $\g\in(0,1)$, we can define the metric 
$$d_\g(\{ x_i \}_{i\ge 1},\{ y_i \}_{i\ge 1})=\g^{k-1}$$
where 
$$k=\inf\{
i\ge 1\st x_i\not=y_i
\}   ,   $$
with the convention that the $\inf$ of the empty set is $+\infty$ and $\g^{+\infty}=0$. 

We define the shift $\s$ as 
$$\fun{\s}{\Sigma}{\Sigma},\qquad
\fun{\sigma}{\{ x_1,x_2,x_2,\dots \}}{\{ x_2,x_3,x_4,\dots \}}  .   $$
If $x_1,\dots,x_l\in(1,\dots,p)$, we define the cylinder 
$$[x_1\dots x_l]=\{
\{ y_i \}_{i\ge 1}\in\Sigma\st y_i=x_i\txt{for}i\in(1,\dots,l)
\}  .  $$
We also set 
$$F_{x_1\dots x_l}=F_{x_1}\circ\dots\circ F_{x_l}  \eqno (1.3)$$
and
$$[x_1\dots x_l]_{\tilde G}=F_{x_1}\circ F_{x_{2}}\circ\dots\circ F_{x_l}(\tilde G)  .  \eqno (1.4)$$
If $x=(x_1x_2\dots)$ and $i\in(1,\dots,p)$ we set $(ix)=(ix_1x_2\dots)$. Now (1.4) implies that 
$$F_i([x_1\dots x_l]_{\tilde G})=[i x_1\dots x_l]_{\tilde G}  .  \eqno (1.5)$$
Since the maps $F_i$ are continuous and $\tilde G$ is compact, the sets 
$[x_1\dots x_l]_{\tilde G}$ are compact. By (1.2) we have that 
$F_i(\tilde G)\subset\tilde G$ for $i\in(1,\dots,p)$; by (1.4) this implies that, for all 
$\{ x_i \}_{i\ge 1}\in\Sigma$, 
$$[x_1\dots x_{l-1} x_l]_{\tilde G}\subset [x_1\dots x_{l-1}]_{\tilde G}
\subset\tilde G  .  $$
From (1.1), (1.3) and (1.4) we get that 
$${\rm diam} ([x_1\dots x_l]_{\tilde G})\le\eta^{l}\cdot
{\rm diam}(\tilde G)   .   \eqno (1.6)$$
Let $\{ x_i \}_{i\ge 1}\subset\Sigma$; by the last two formulas and the finite intersection property we have that 
$$\bigcap_{l\ge 1}[x_1\dots x_l]_{\tilde G}$$
is a single point, which we call $\Phi(\{ x_i \}_{i\ge 1})$; formula (1.6) implies in a standard way that the map $\fun{\Phi}{\Sigma}{\tilde G}$ is continuous. It is not hard to prove, using (1.2), that $\Phi$ is surjective. We shall call $\tilde d$ the distance on $\tilde G$ induced by the Euclidean distance on $\R^d$ and, from now on, in our choice of the metric on 
$\Sigma$ we take $\g\in(\eta,1)$; this implies by the definition of $d_\g$ on $\Sigma$ and (1.6) that $\Phi$ is 1-Lipschitz.  

In [4] we needed some control on the lack of injectivity of $\Phi$; since we are going to apply this theory to the limit set $\tilde G$ of the sequence (5), which is totally disconnected, injectivity comes for free and we can introduce some simplifications 
{\it vis-\`a-vis} [4]. 

\noindent {\bf (F2)} We ask that, for $i\not=j$, $F_i(\tilde G)\cap F_j(\tilde G)$ is empty. 

It is easy to see that (F2) implies that $\tilde G$ is totally disconnected and the map 
$\Phi$ is injective.

\noindent {\bf (F3)} Since by (F2) the sets $F_i(\tilde G)$ are disjoint and compact and since (1.2) holds,  we can find disjoint open sets $\oc_1,\dots,\oc_p\subset\R^d$ such that
$$\tilde G\cap\oc_i=F_i(\tilde G) \txt{for}i\in(1,\dots,p)  .  $$
We ask that $\oc_i\subset F_i^{-1}(\oc_i)$ (or, equivalently, that 
$F_i(\oc_i)\subset\oc_i$, since the maps $F_i$ are diffeos). 

We define a map $\fun{F}{\bigcup_{i=1}^p\oc_i}{\R^d}$ by
$$F(x)=F_i^{-1}(x)  \txt{if}x\in\oc_i  .  $$
This implies the second equality below; the first one comes since we supposed that 
$\oc_i\subset F_i^{-1}(\oc_i)$ for all $i\in(1,\dots,p)$. 
$$F\circ F_i(x)=x\qquad\forall x\in\oc_i\subset F_i^{-1}(\oc_i)
\txt{and}
F_i\circ F(x)=x\qquad\forall x\in\oc_i  .  \eqno (1.7)$$
We call $a_i$ the unique fixed point of $F_i$; note that, by (1.2), $a_i\in G$. 

\vskip 1pc

Note that, if $x=(x_1x_2\dots)$, then by the definitions of $\Phi$ and $\s$ 
$$\Phi(x)=\bigcap_{l\ge 1}[x_1\dots x_l]_{\tilde G},\qquad
\Phi\circ\s(x)=\bigcap_{l\ge 2}[x_2\dots x_l]_{\tilde G}  .  $$
This implies the first and last equalities below. Now recall that, if $x=(x_1x_2\dots)$, then 
$\Phi(x)\in F_{x_1}(\tilde G)$ by (F2); thus, $\Phi(x)\in\oc_{x_1}$ and (1.4) and (1.7) imply the middle equality. 
$$F\circ\Phi(x)=F\left(
\bigcap_{l\ge 1}[x_1\dots x_l]_{\tilde G}
\right)  =    
\bigcap_{l\ge 2}[x_2\dots x_l]_{\tilde G}=
\Phi\circ\s(x)  .  $$
As a consequence, the two equalities below hold for all $x\in\Sigma$. 
$$\left\{\matrix{
\Phi\circ\s(x)=F\circ\Phi(x)\cr
\Phi(i,x) =F_i(\Phi(x))   
\quad\forall i\in (1,\dots,p)   .   
}
\right.       \eqno (1.8)  $$
In other words, up to a change of coordinates, shifting the coding one place to the left is the same as applying $F$. Iterating the first one of (1.8) we get that, for all 
$l\ge 1$ and all 
$x\in\Sigma$, 
$$\Phi\circ\s^l(x)=F^l\circ\Phi(x)   .    $$

A particular case (save for the fact that (F2) is not satisfied, but the more general hypotheses of [4] are) is the harmonic Sierpinski gasket on $\R^2$; we refer the reader to [9] for an introduction to the properties of this set. We set
$$T_1=\left(
\matrix{
\frac{3}{5},&0\cr
0,&\frac{1}{5}
}
\right)  ,  \quad
T_2=\left(
\matrix{
\frac{3}{10},&\frac{\sqr 3}{10}\cr
\frac{\sqr 3}{10},&\frac{1}{2}
}
\right)  ,\quad
T_3=\left(
\matrix{
\frac{3}{10},&-\frac{\sqr 3}{10}\cr
-\frac{\sqr 3}{10},&\frac{1}{2}
}
\right)   ,  $$
$$A=\left(
\matrix{
0\cr
0
}
\right),\qquad
B=\left(
\matrix{
1\cr
\frac{1}{\sqrt 3}
}
\right) ,\qquad
C=\left(
\matrix{
1\cr
-\frac{1}{\sqrt 3}
}
\right)  $$
and 
$$F_1(x)=T_1(x),\quad F_2(x)=
B   +T_2\left(
x-B
\right),\quad
F_3(x)=
C   +T_3\left(
x-C
\right)  .  $$
Referring to figure 2 below, $F_1$ brings the triangle $ABC$ into $Abc$; 
$F_2$ brings $ABC$ into $Bac$ and $F_3$ brings $ABC$ into $Cba$. We take three disjoint open sets $\oc_1$, $\oc_2$, $\oc_3$ which contain, respectively, the triangle 
$Abc$ minus $b,c$, $Bca$ minus $c,a$ and $Cba$ minus $a,b$. 

We define the map $F$ as 
$$F(x)=F_i^{-1}(x)
\txt{if}
x\in\oc_i $$
and on $\{ a,b,c \}$ we extend it arbitrarily, say $F(a)=A$, $F(b)=B$ and $F(c)=C$.

\centerline{Figure 2}

\vskip 1pc

We have stated these hypotheses essentially with theorem 1.1. below in mind; before coming to it, we need some notation on matrix-valued measures. This is because, recently, several authors ([8], [13], [17]) realised that Kusuoka's measure $\kappa$ ([12]) is an object well-known in Dynamical Systems, namely a Gibbs measure; we shall follow the approach of [4] and [5] in which the carr\'e\ de champs is expressed  through a matrix-valued Gibbs measure $\tau$.

\noindent{\bf The space of matrices and the measures valued in it.} We define $\tilde M$ as the space of all matrices from $\R^2$ to itself; we call $M$ the subspace of symmetric matrices. The space $\tilde M$ is a Hilbert space under the Hilbert-Schmidt inner product 
$$(A,B)_{HS}\colon={\rm tr}(A\tr B)$$
where $\tr B$ denotes the transpose of $B$. The norm is the standard one:
$$||A||_{HS}^2\colon=(A,A)_{HS}   .   $$

We say that $A\in M$ is positive (or semi-positive) definite if $(Av,v)>0$ (or $(Av,v)\ge 0$) for all $v\in\R^2\setminus\{ 0 \}$. 

For positive or semipositive symmetric matrices we shall use the standard notation, i. e. $A>0$ and $A\ge 0$ respectively. If $A,B\in M$, we shall say that $B\ge A$ if 
$B-A\ge 0$. 

Let $G\subset\R^2$ be compact; we define $\mc(G,\tilde M)$ as the space of the Borel measures on $G$ valued in $\tilde M$. We define the integral of $A\in C(G,\tilde M)$ against 
$\mu\in\mc(G,\tilde M)$; in order to do this, we recall from [18] that the total variation of 
$\mu$ is a finite scalar measure 
$||\mu||$ on the Borel sets of $G$. The polar decomposition of $\mu$ is given by 
$$\mu=M_x||\mu||$$
where $||\mu||$ is the total variation measure of $\mu$ and 
$\fun{M}{G}{\tilde M}$ is a Borel field of matrices which satisfies 
$$||M_x||_{HS}=1
\txt{for $||\mu||$-a. e. $x\in G$.}  \eqno (1.9)$$
If $\fun{A}{G}{\tilde M}$ is Borel and $||A||_{HS}\in L^1(G,||\mu||)$, then by (1.9) and Cauchy-Schwarz we have that 

\noindent $(A_x,M_x)_{HS}\in L^1(G,||\mu||)$. Consequently, we can define the scalar 
$$\int_G(A_x,\dr\mu(x))_{HS}\colon =
\int_G(A_x,M_x)_{HS}\dr||\mu||(x)  .  $$
The duality coupling between $C(G,\tilde M)$ and $\mc(G,\tilde M)$ is given by 
$$\fun{\inn{\cdot}{\cdot}}{
C(G,\tilde M)\times\mc(G,\tilde M)
}{\R}  $$
$$\inn{A}{\mu}\colon =
\int_G(A_x,\dr\mu(x))_{HS}  .  \eqno (1.10)$$
If $Q\in C(G,\tilde M)$ and $\mu\in\mc(G,\tilde M))$, we define the scalar measure 
$(Q,\mu)_{HS}$ by 
$$\int_G f(x)\dr(Q,\mu)_{HS}(x)\colon =
\int_G(fQ,\dr\mu)_{HS}    \eqno (1.11)$$
for all $f\in C(G,\R)$. In other words, $(Q,\mu)_{HS}=(Q_x,M_x)_{HS}\cdot ||\mu||$ for the decomposition (1.9). 

If $\fun{v,w}{G}{\R^d}$ are Borel functions such that 
$||v||\cdot ||w||\in L^1(G,||\mu||)$ then, again by (1.9) and Cauchy-Schwarz, 
$(v_x,M_xw_x)\in L^1(G,||\mu||)$, which implies that the third integral below converges. The second equality below is the definition of the middle integral, the first one comes from the fact that $(v,Mw)=(v\otimes w,M)_{HS}$.  
$$\int_G(v_x\otimes w_x,\dr\mu(x))_{HS} =
\int_G(v_x,\dr\mu(x)\cdot w(x))=
\int_G(v_x,M_x w_x)\dr||\mu||(x)  .  \eqno (1.12)$$

Now we concentrate on symmetric matrices, i. e. on the space $M$ defined above; the reason is that Riemannian tensors are symmetric and we want to define a Riemannian structure on $G$, natural for the dynamics $F$ of (1.7).  

We define $\mc(G,M)$ as the linear space of the Borel measures on $G$ valued in $M$. Since $M\subset\tilde M$ we have that  
$$\mc(G,M)\subset\mc(G,\tilde M)   .  \eqno (1.13)$$
By Lusin's theorem, $\mu\in\mc(G,\tilde M)$ belongs to $\mc(G,M)$ if and only if, for all 
$\a,\b\in C(G,\R^2)$ and the duality coupling (1.10) we have that 
$$\langle\a\otimes\b,\mu\rangle=\langle\b\otimes\a,\mu\rangle  .  \eqno (1.14) $$

We say that $\mu\in\mc(G,M)$ is semipositive definite if $\mu(E)$ is semipositive-definite  for all Borel sets $E\subset S$. Using again Lusin's theorem we see as in [4] that 
$\mu\in\mc(G, M)$ is semipositive definite if and only if 
$$\inn{A}{\mu}\ge 0   \eqno (1.15)$$
for all $A\in C(G,M)$ such that $A_x\ge 0$ for all $x\in G$. 

We denote by $\mc_+(G,M)$ the set of all semipositive definite measures of 
$\mc(G,M)$; by the characterisations (1.14)-(1.15), $\mc_+(G,M)$ is a convex set of $\mc(G,\tilde M)$, closed for the  weak$\ast$ topology. 

Let now $Q\in C(G,M)$ be such that $Q_x$ is positive-definite for all $x\in G$; since $G$ is compact there is $\e>0$ such that 
$$\e Id\le Q_x\le\frac{1}{\e}Id\qquad\forall x\in G  .  \eqno (1.16)$$
For such a matrix $Q$ we define $\pc_Q(G,M)$ as the set of all 
$\mu\in\mc_+(G,M)$ such that 
$$\int_G(Q_x,\dr\mu(x))_{HS}=1  .  $$
As shown in [4], if $Q$ satisfies (1.16) there is $D_1(\e)>0$ such that, for all 
$\mu\in\mc_+(G,M)$, 
$$\frac{1}{D_1(\e)}(Q,\mu)_{HS}\le
||\mu||\le D_1(\e)(Q,\mu)_{HS}   .  \eqno (1.17)$$
As a consequence, the two measures $||\mu||$ and $(Q,\mu)_{HS}$ are mutually absolutely continuous and their mutual Radon-Nikodym derivatives are bounded. Integrating (1.17) over $G$ we get that, if $\mu\in\pc_Q(G,M)$, then
$$||\mu||(G)\le D_1(\e)  .  $$
By its definition and formulas (1.14)-(1.15), $\pc_Q(G, M)$ is a convex subset of 
$\mc(G,\tilde M)$, closed for the weak$\ast$ topology; by the formula above, it is compact. 

We introduce a last bit of notation. Let the set $\tilde G$ be as in (1.2); we fix a compact set $G_0$ and we define iteratively the pre-fractal $G_l$ as 
$$G_{l}=\bigcup_{i=1}^p F_i(G_{l-1})
\txt{for}l\ge 1  .  $$
By [10], $G_l$ converges to the fractal $\tilde G$ in the Haussdorff distance for sets. 

In particular, we can fix $R>0$ such that $B(0,R)$ contains all the pre-fractals $G_l$; since the maps $F_i$ are contractions, possibly enlarging $R$ we can require that 
$F_i(\bar B(0,R))\subset\bar B(0,R)$ for all $i\in(1,\dots,p)$. Let 
$Q\in C(\bar B(0,R),M)$ which satisfies (1.16) on $\bar B(0,R)$; we define 
$\pc_Q(\bar B(0,R),M)$ as the set of the Borel measures on $\bar B(0,R)$ which are valued in $M$, are positive-definite in the sense above and satisfy 
$$\int_{B(0,R)}(Q_x,\dr\mu(x))_{HS}=1  .  $$ 

\noindent{\bf The Ruelle operator.} Let $\{ F_i \}_{i=1}^p$ be the affine maps defining the fractal and let $V\in C^{0,\alpha}(\tilde G,\R)$ for some $\alpha\in(0,1]$. We define the Ruelle operator 
$$\fun{\L}{C(\tilde G,M)}{C(\tilde G,M)}$$ 
by
$$(\L A)(x)\colon =\sum_{i=1}^p e^{V\circ F_i(x)}\cdot
\tr(DF_i)\cdot A_{F_i(x)}\cdot (DF_i)    .  \eqno (1.18)$$
Since $F_i$ is affine, its derivative is constant and we haven't marked the point where we calculate it. 

It is easy to see that $\L$ is continuous; its adjoint 
$$\fun{\L^\ast}{\mc(\tilde G,M)}{\mc(\tilde G,M)} $$
is defined by 
$$\inn{\L A}{\mu}=\inn{A}{\L^\ast\mu}  $$
for all $A\in C(\tilde G,\tilde M)$ and $\mu\in\mc(\tilde G,\tilde M)$; the duality coupling is that of (1.10). 

In Dynamical Systems theory, the standard procedure to find Gibbs measures is to apply the Perron-Frobenius theorem to $\L$ and a suitable cone contained in $C(\tilde G,M)$; in [4], [5] it is shown that this is possible if the following nondegeneracy hypothesis holds.

\vskip 1pc

\noindent {\bf (ND)} We suppose that there is $\g>0$ such that for all 
$c,e\in \R^d$ we can find 
$i\in(1,\dots,p)$ such that 
$$|((DF_i) c, e)|\ge\g ||c||\cdot||e||    $$
where the inner product on the left is that of $\R^d$. 

In other words we assume that, given any $c\in\R^d\setminus\{ 0 \}$, the vectors 
$\{ (DF_i) c \}_{1\le i\le p}$ generate $\R^d$. 

\vskip 1pc

It is immediate that $\L$ preserves the cone of semi-positive definite matrices. Actually, it is easy to see that $\L$ sends this cone in a cone strictly contained in it: namely, (ND) implies that, if $A$ is a field of semi-positive definite matrices, then $\L A$ is a field of positive-definite matrices. Now it is easy to apply ([4], [5]) the Perron-Frobenius theorem as in [16] and [20] to show the following theorem; we state the version of $C(B(0,R),M)$ instead of that on $C(\tilde G,M)$. 

\thm{1.1} Let us suppose that (F1)-(F3) and (ND) hold; let $V\in C^{0,\a}(B(0,R),\R)$ with 
$\a\in(0,1]$. Then, we have the following. 

\noindent 1) The operator $\L$ of (1.18) has a simple, positive eigenvalue $\l>0$ with a positive-definite eigenvector. In other words, there is 
$Q\in C(\bar B(0,R),M)$ such that $Q(x)$ is positive-definite for all $x\in\bar B(0,R)$ and 
$$\L Q=\l Q  .  $$
Moreover, if $\tilde Q\in C(\bar B(0,R),M)$ satisfies $\L\tilde Q=\l\tilde Q$, then 
$\tilde Q=\d Q$ for some $\d\in\R$. 

\noindent 2) The number $\l$ is an eigenvalue also for $\L^\ast$. More precisely, let us fix a positive eigenvector $Q$ as in point 1) and let $\pc(\bar B(0,R),M)$ be the set we defined after (1.17). Then, the map $\frac{1}{\l}\L^\ast$ brings $\pc(\bar B(0,R),M)$ into itself and there is a unique $\tau\in\pc_{Q}(\bar B(0,R),M)$ such that 
$$\L^\ast\tau=\l\tau  .  $$
The measure $\tau$ is supported on $\tilde G$. 

\noindent 3) Let us define the scalar measure $\kappa\colon=(Q,\tau)_{HS}$, where the notation is that of (1.11). Then, $\kappa$ is a probability measure on $\tilde G$ ergodic for the expansive map $F$ defined in (F3). Moreover, $\kappa$ is non-atomic; $\kappa$ and 
$\tau$ are mutually absolutely continuous, with bounded Radon-Nikodym derivatives 
$\frac{\dr\kappa}{\dr||\tau||}$ and $\frac{\dr\tau}{\dr\kappa}$. 

\noindent 4) By (1.13), for $f,g\in C(\R^d,\R^d)$ we can define 
$$\ec^1(f,g)\colon=\int_{\tilde G}(( f(x)\otimes  g(x),\dr\tau(x))_{HS}=
\int_{\tilde G}(f(x),\dr\tau\cdot g(x))    \eqno (1.19)  $$
where the integrals are as in (1.12). If $V\equiv 0$, then $\ec^1$ is self-similar; in other words, for the maps $F_i$ of (1.1) and all $f,g\in C^1(\R^d,\R^d)$ we have that 
$$\ec^1(f,g)=\frac{1}{\l}\sum_{i=1}^d\ec^1((DF_i) f,(DF_i) g)   .   \eqno (1.20)$$

\noindent 5) Let $\tilde\tau$ be any element of $\pc_Q(\bar B(0,R),M)$ and let $\tau$ be the eigenvector of point 2); then,  
$$\left(
\frac{1}{\l}\L^\ast
\right)^l\tilde\tau\tends \tau \txt{for}l\tends+\infty$$
in the weak$\ast$ topology of $\mc(G,M)$. 

\noindent 6) The measure $\tau$ has the Gibbs property; in other words, with the notation of (1.3) and (1.4), there is a positive constant $D_1>0$ such that, for all $l\ge 1$ and all 
$x=\{ x_i \}_{i\ge 0}\in\Sigma$, 
$$\frac{e^{V_l(\Phi(x))}}{D_1\l^l}(DF_{x_1\dots x_{l}})\cdot
\tau(\tilde G)\cdot\tr(DF_{x_1\dots x_{l}})\le
\tau([x_1\dots x_{l}]_{\tilde G})\le$$
$$\frac{D_1 e^{V_l(\Phi(x))}}{\l ^l}
(DF_{x_1\dots x_{l}})\cdot
\tau(\tilde G)\cdot\tr(DF_{x_1\dots x_{l}})     $$
where 
$$V_l(\Phi(x))\colon= 
V(\Phi(x))+V(\Phi\circ\s(x))+\dots+V(\Phi\circ\s^{l-1}(x))    $$
and $F_{x_1\dots x_{l}}$ is as in (1.3).

\noindent If $V\equiv 0$, then $D_1=1$ and the inequalities above are equalities; in the case of the harmonic Sierpinski gasket, $\tau(G)=Id$.

\rm

\vskip 1pc

Another fact, also shown in [4], is that there is a parallel Ruelle operator on the coding space; namely, for the space $\Sigma$ we defined after (1.2) we can set 
$$\fun{\L_\Sigma}{C(\Sigma,\R)}{C(\Sigma,\R)}$$
$$(\L_\Sigma A)(x)=\sum_{i=1}^p
e^{V\circ\Phi(ix)}\cdot\tr(DF_i)\cdot A_{(ix)}\cdot (DF_i)  .   \eqno (1.21)$$
Note that $V\circ\Phi\in C^{0,\a}(\Sigma,\R)$; indeed, we saw above that that the coding map $\Phi$ is Lipschitz and we are supposing that $V\in C^{0,\a}(G,\R)$. 

For $\L_\Sigma$ there is an analogue of theorem 1.1, and the relationship between the eigenvalues and eigenvectors of $\L$ on $C(G,M)$ and $\L_\Sigma$ on $C(\Sigma,M)$ is the natural one: for instance ([4]), the Perron-Frobenius eigenvalue of $\L_\Sigma$ coincides with that of $\L$. Moreover, if we multiply the eigenvector $Q_\Sigma$ of 
$\L_\Sigma$ by a suitable constant we get that $Q_\Sigma=Q\circ\Phi$.  Lastly, the eigenvector $\tau_\Sigma$ of $\L_\Sigma^\ast$ and Kusuoka's measure $\kappa_\Sigma$ on $\Sigma$ satisfy 
$$\tau=\Phi_\sharp\tau_\Sigma,\qquad \kappa=\Phi_\sharp\kappa_\Sigma      $$
where $\Phi_\sharp$ denotes the push-forward of measures by $\Phi$. 

In the following, we shall need a variation of this construction. Namely, let 
$$\{
F^i_1,\dots,F^i_p
\}_{i\ge 1}   $$
be a sequence of $p$-uples of affine contractions of $\R^d$; we suppose that they are all 
$\eta$-Lipschitz for some common $\eta\in(0,1)$ and that there is a compact set 
$\hat G_0\subset\R^d$ with the following two properties: 

\noindent 1) $F^i_j(\hat G_0)\subset\hat G_0$ for all $i\ge 1$ and $j\in(1,\dots,p)$, and 

\noindent 2) For all $i\ge 1$, $F^i_j(\hat G_0)\cap F^i_l(\hat G_0)=\emptyset$ if $j\not=l$. 

We define 
$$F_{i_1\dots i_{l}}\colon=
F_{i_1}^1\circ F^2_{i_2}\circ\dots\circ F^{l}_{i_{l}}  ,  \eqno (1.22)$$
and
$$[i_1\dots i_{l}]_{\hat G_0}\colon=
F_{i_1\dots i_{l}}(\hat G_0)   .    $$
Let $i=(i_1,i_2,\dots)\in\Sigma$; point 1) above implies that  
$$[i_1\dots i_{l}]_{\hat G_0}\subset [i_1\dots i_{l-1}]_{\hat G_0}$$ 
and thus the finite intersection property holds. Together with the fact that the diameter of 
$[i_1\dots i_l]_{\hat G_0}$ tends to zero, we get that there is a unique point 
$\Phi(i)\in\R^2$ such that  
$$\{ \Phi(i) \}=
\bigcap_{i\ge 1}[i_1\dots i_{l}]_{\hat G_0}   .    \eqno (1.23)$$
It is easy to see that $\Phi(i)$ does not depend on the choice of the set $\hat G_0$, provided this set satisfies points 1) and 2) above. 

We skip the easy proof that 2) implies that $\Phi$ is a homeomorphism of $\Sigma$ onto its image. Thus, if we set 
$$\tilde G=\Phi(\Sigma)     \eqno (1.24)$$
we get that $\tilde G$ is totally disconnected.

The Ruelle operator behaves well in this situation if we suppose, as we shall do in the rest of the paper, that $V\equiv 0$ and that the $p$-uples $(F^i_1,\dots, F^i_p)$ are homotheties of the first one, $(F^1_1,\dots, F^1_p)$. More precisely, we suppose that 
$$DF^i_j=\e_i DF^1_i
\txt{for all}i\ge 1,\quad j\in(1,\dots,p)    $$
and we call $\L_i$ the Ruelle operator on $C(\Sigma,M)$ defined as in (1.21), but for the maps $\{ F^i_j \}_{j\in(1,\dots,p)}$. Now (1.21) immediately implies that 
$\L_i=\e_i^2\L_1$; thus, the eigenvalue $\l_i$ of $\L_i$ satisfies $\l_i=\e_i^2\l_1$. As a consequence, all the operators $\frac{1}{\l_i}\L_i$ on $C(\Sigma,\R)$ coincide with 
$\frac{1}{\l_1}\L_1$. As we recalled after theorem 1.1, $\frac{1}{\l_1}\L_1$ induces a measure-valued Gibbs measure $\tau_\Sigma$ on $\Sigma$, together with a scalar Kusuoka measure 
$\kappa_\Sigma$. We call $\Phi_\rc$ the conjugation of (1.23) associated to the sequence $\rc=\{ \e_i \}_{i\ge 1}$ (we are abusing the notation of (7)) and we define two measures on $\tilde G$: 
$$\tau_\rc\colon=(\Phi_\rc)_\sharp\tau_\Sigma,\qquad 
\kappa_\rc\colon=(\Phi_\rc)_\sharp\kappa_\Sigma   .   \eqno (1.25)$$
It follows easily by the definition of push-forward that these two measures are supported on the set $\tilde G$ of (1.24). 

We define 
$$\fun{\tilde F}{\Phi_\rc(\Sigma)}{\Phi_{\s\rc}(\Sigma)}$$
$$\tilde F(x)=(F^1_i)^{-1}(x)
\txt{if}x\in F^0_i(\tilde G)   .   $$
We can see as in (1.8) that $\Phi_{\s\rc}\circ\s=\tilde F\circ\Phi_{\rc}$. Since $\kappa_\Sigma$ is $\s$-invariant (i. e. $\s_\sharp\kappa_\Sigma=\kappa_\Sigma$), formula (1.25) implies that 
$$\tilde F_\sharp\kappa_{\rc}=\kappa_{\s\rc}    .   $$
We define the energy as in (1.19), i. e. 
$$\ec_\rc^1(\nabla u,\nabla v)=
\int_{\tilde G}(\nabla u,\dr\tau_\rc\cdot\nabla v)   .  $$

Let now $G_0\subset\R^d$ be a compact set, let $G_l$ be defined by 
$$G_l=\bigcup_{i_1\dots i_l\in\{ 1,\dots,p \}}^p F_{i_1\dots i_l}(G_{0})  $$
and let $B(0,R)$ be a ball which contains all the sets $G_l$. Let us set 
$$\fun{\tilde\L_{\tilde G,l}}{C(\bar B(0,R),M)}{C(\bar B(0,R),M)}$$
$$(\tilde\L_{\tilde G,l} A)(x)=\sum_{i=1}^p 
\tr(DF^l_i)\cdot A_{F^l_i(x)}(DF^l_i)  .  \eqno (1.26)$$
We can restrict this operator to $C(\tilde G,M)$; since 
$\Phi_\rc(ix)=F^1_i\circ\Phi_{\s\rc}(x)$, there is a relationship between the Ruelle operators on $\Sigma$ and $\tilde G$:  
$$\L_1(A\circ\Phi_{\rc})=(\L_{\tilde G_,1}A)\circ\Phi_{\s\rc}   .   $$
Since $\tau_\Sigma$ is an eigenvector for $\L_1^\ast$, the last formula implies that 
$$(\Phi_{\rc})_\sharp\tau_\Sigma=
\frac{1}{\l_1}(\L_{\tilde G,1})^\ast(\Phi_{\s\rc})_\sharp\tau_\Sigma  .  $$
We define 
$$\ec^1_\rc(f,g)=\int_G(f,\dr\tau_\rc\cdot g)   .  $$
The last two formulas imply that, if $u,v\in C^1(\R^d,M)$, then 
$$\ec_\rc^1(\nabla u,\nabla v)=
\frac{1}{\l_1}\sum_{i=1}^p
\ec_{\s\rc}(\nabla(u\circ F^1_i),\nabla(v\circ F^1_i))   .   \eqno (1.27)$$
Let $\tilde\tau_0\in\pc_Q(\Sigma,M)$; since the operators $\frac{1}{\l_i}\L_i$ coincide, theorem 1.1 implies that 
$$\left( \frac{1}{\l_1}\L_1^\ast \right)\cdot\dots \cdot 
\left( \frac{1}{\l_l}\L_l^\ast \right)\tilde\tau_0\tends\tau_\Sigma  .  $$
Let $\tau_0\in\pc_Q(\bar B(0,R),M)$ be supported on $G_0$; we set 
$$\tau_l=\left(
\frac{1}{\l_1}\L_{\tilde G,1}^\ast
\right)\cdot\dots\cdot
\left(
\frac{1}{\l_l}\L_{\tilde G,l}^\ast
\right)  \tau_{0}  .  \eqno (1.28)$$
It is easy to see, using the definition of $\L_{\tilde G,l}$, that $\tau_l\in \pc_Q(B(0,R),M)$ and that $\tau_l$ is supported on the set $G_l$ defined before (1.26). Moreover, each operator $\frac{1}{\l_l}\L_{\tilde G,l}$ is a contraction of the cone ${\cal C}$ of positive-definite matrices in the hyperbolic metric, its Lipschitz constant is bounded away from 1 and all sets $\frac{1}{\l_l}\L_{\tilde G,l}({\cal C})$ are contained in the same bounded subcone of ${\cal C}$. Using these facts it is easy to prove that, with the notation we used above, 
$$\tau_l\tends\tau_\rc$$
for the measure  $\tau_\rc$ defined above. 

In particular, defining  $\ec_\rc(\nabla u,\nabla v)$ as above, and setting 
$$\ec_{\rc,l}^1(\nabla u,\nabla v)=
\int_{\bar B(0,R)}(\nabla u,\dr\tau_l\cdot\nabla v)$$
we get that, for all $u,v\in C^1(\R^2,\R)$ we have that 
$$\ec_{\rc,l}^1(\nabla u,\nabla v)\tends\ec_\rc(\nabla u,\nabla v)  .  \eqno (1.29)$$

\vskip 2pc

\centerline{\bf \S 2}
\centerline{\bf The affine maps}

\vskip 1pc

We could prove theorem 1 verifying that the maps $F^i_1$ of (3) (together with $F^i_2$ and $F^i_3$ defined by symmetry) generate harmonic pre-fractals. We follow a different road, considering maps $F_1^i$ more general than the ones of (3). In the next two sections we shall see that, if we want the pre-fractals to be harmonic, then the maps $F^i_1$ must have the form (3).

\noindent{\bf Definition by symmetry.} Let $\a,\b\in(0,1)$; using the column notation for vectors we set 
$$T_1\left(
\matrix{
x\cr
y
}
\right)  =
\left(
\matrix{
\a x\cr
\b y
}
\right)  $$
$$F_1\left(
\matrix{
x\cr
y
}
\right)  =
T_1\left(
\matrix{
x\cr
y
}
\right)  .  \eqno (2.1)$$
The function $F_1$ is an affine (actually, linear) contraction which fixes point $A$ of (1). Now we define by symmetry the other two maps $F_2$ and $F_3$, which fix $B$ and $C$ respectively. 

We denote by $R_\theta$ the rotation of angle $\theta$ in the anticlockwise direction, i. e. 
$$R_\theta=\left(
\matrix{
\cos\theta, &-\sin\theta\cr
\sin\theta, &\cos\theta
}
\right)  .  \eqno (2.2)$$
We set 
$$T_2\left(
\matrix{
x\cr
y
}
\right)=
R_{\frac{-2\pi}{3}}\cdot T_1\cdot R_{\frac{2\pi}{3}}\left(
\matrix{
x\cr
y
}
\right)  ,  $$
$$F_2\left(
\matrix{
x\cr
y
}
\right)=
T_2\left[
\left(
\matrix{
x\cr
y
}
\right)-B
\right] +B,   \eqno (2.3)$$
$$T_3\left(
\matrix{
x\cr
y
}
\right)=
R_{\frac{2\pi}{3}}\cdot T_1\cdot R_{-\frac{2\pi}{3}}\left(
\matrix{
x\cr
y
}
\right)  ,  $$
$$F_3\left(
\matrix{
x\cr
y
}
\right)=
T_3\left[
\left(
\matrix{
x\cr
y
}
\right)-C
\right] +C  .   \eqno (2.4)$$

\noindent{\bf Disconnection.} Let now $\rc=\{ (\a_i,\b_i) \}_{i\ge 1}$ and let 
$F^i_1$, $F_2^i$, $F^i_3$ be the maps of (2.1), (2.3) and (2.4) respectively for the parameters $(\a_i,\b_i)$. We want to find conditions on the parameters such that conditions 1) and 2) before formula (1.22) hold. As we saw there, this implies that the sequence 
$\{ F^i_1,F^i_2,F^i_3 \}_{i\ge 1}$ determines a set $\tilde G$ by (1.24); since for all $i$ 
$F^i_1$ fixes $A$, $F^i_2$ fixes $B$ and $F^i_3$ fixes $C$, it is easy to see that 
$A,B,C\in\tilde G$ and that $F_{i_1\dots i_l}(j)\in\tilde G$ for all $j\in\{ A,B,C \}$ and 
$i_1,\dots,i_l\in\{ 1,2,3 \}$. Moreover, $\tilde G$ is totally disconnected since, by (1.24), it is homeomorphic to the totally disconnected set $\Sigma$. 

Let us denote by $\hat G_0$ the compact set bounded by the triangle $G_0$ of (4) (the solid triangle, if you want). If $\a_i\ge\b_i$ and $\a_i,\b_i\in(0,1)$ then the map $F_1^i$ brings $\hat G_0$ into itself; since the maps $F_2^i$ and $F_3^i$ share the same property by symmetry, point 1) follows. Note that, by (2.1), 
$F_1^i(\hat G_0)\subset\max(\a_i,\b_i)\hat G_0$ and a similar inclusion holds for $F^i_2(\hat G_0)$; thus,  $F^i_1(\hat G_0)\cap\hat F^i_2(\hat G_0)=\emptyset$ if 
$$[B+\max(\a_i,\b_i)(\hat G_0-B)]\cap\max(\a_i,\b_i)\hat G_0=\emptyset   .    $$
Since $\hat G_0$ is an equilateral triangle with side 1, this is implied by 
$ \a_i,\b_i\in\left( 0,\2\right)$. Summing up, points 1) and 2) are implied by the two conditions below. 
$$\a_i\ge\b_i,\qquad \a_i,\b_i\in\left( 0,\2\right)
\txt{for all}i\ge 0  .  \eqno (2.5)$$
It will be possible to relax the second condition above for the particular choice of 
$(\a_i,\b_i)$ we shall make in section 4 below. 

\noindent{\bf Explicit formula for $F_2^i$.}  We begin with an explicit expression for 
$T_2^i$. The first equality below is the formula before (2.3), the second one follows from (2.2) with $\theta=\frac{2\pi}{3}$. 
$$T_2^i=R_{-\frac{2\pi}{3}}\cdot T_1^i\cdot R_{\frac{2\pi}{3}}=$$
$$\left(
\matrix{
-\2, &\frac{\sqrt 3}{2}\cr
\frac{-\sqrt 3}{2},&-\2
}
\right)
\left(
\matrix{
\a_i,&0\cr
0,&\b_i
}
\right)
\left(
\matrix{
-\2, &\frac{-\sqrt 3}{2}\cr
\frac{\sqrt 3}{2},&-\2
}
\right)  =$$
$$\left(
\matrix{
\frac{\a_i+3\b_i}{4},&\frac{\sqrt 3}{4}(\a_i-\b_i)\cr
\frac{\sqrt 3}{4}(\a_i-\b_i),&\frac{3\a_i+\b_i}{4}
}
\right)  .  $$
The first equality below follows from (2.3) and (1), the third one comes from the formula above and again (1). 
$$F_2^i(A)=T_2^i(-B)+B=(Id-T_2^i)(B)=$$
$$\left(
\matrix{
\frac{4-(3\b_i+\a_i)}{4},&\frac{\sqrt 3}{4}(\b_i-\a_i)\cr
\frac{\sqrt 3}{4}(\b_i-\a_i),&\frac{4-(\b_i+3\a_i)}{4}
}
\right)
\left(
\matrix{
\frac{\sqrt 3}{2}\cr
\2
}
\right)  =
\left(
\matrix{
\frac{\sqrt 3}{4}(2-\b_i-\a_i)\cr
\frac{1}{4}(2+\b_i-3\a_i)
}
\right)   .  \eqno (2.6)$$
By (2.1) and (1) we get that 
$$F_1^i(B)=\left(
\matrix{
\a_i\frac{\sqrt 3}{2}\cr
\frac{\b_i}{2}
}
\right)  .  \eqno (2.7)$$

\noindent{\bf Explicit parametrisation of the edges.} In both formulas below, (8) implies the first equality, (1) the second one.  
$$\g_{BA}(t)=(1-t)B=
(1-t)\left(
\matrix{
\frac{\sqrt 3}{2}\cr
\2
}
\right)  .  \eqno (2.8)$$
$$\g_{BC}(t)=(1-t)B+tC=
\left(
\matrix{
\frac{\sqrt 3}{2}\cr
\2-t
}
\right)  .   \eqno (2.9)$$
We write down explicitly $\g_{F_1^i(B)F_2^i(A)}$; the first equality below comes from (8), the second one from (2.6) and (2.7). 
$$\g_{F_1^i(B)F_2^i(A)}(t)=
(1-t)F_1^i(B)+tF_2^i(A)=$$
$$(1-t)\left(
\matrix{
\a_i\frac{\sqrt 3}{2}\cr
\frac{\b_i}{2}
}
\right)+
t\left(
\matrix{
\frac{\sqrt 3}{4}(2-\b_i-\a_i)\cr
\frac{1}{4}(2+\b_i-3\a_i)
}
\right) =$$
$$\left(
\matrix{
\frac{\a_i\sqrt 3}{2}+t\frac{\sqrt 3}{4}(2-\b_i-3\a_i)\cr
\frac{\b_i}{2}+\frac{t}{4}(2-\b_i-3\a_i)
}
\right)  .  \eqno (2.10)$$

\vskip 2pc\centerline{\bf \S 3}
\centerline{\bf Definition of the energy on the $l$-th pre-fractal }

\vskip 1pc

Before defining harmonicity on the $l$-th pre-fractal we must define the Dirichlet form, which we do using the affine maps of the last section; since at this point the only condition on the parameters $(\a_i,\b_i)$ is (2.5), we are going to get a form more general than that of (9), (10) and (11). In this section we write down the form without explaining one of its features, i. e. its connection with the dynamics induced by the maps $F^i_j$; this connection will be studied in section 5. 

Let us consider a sequence $\{ (\a_i,\b_i) \}_{i\ge 1}\subset(0,+\infty)^2$. For each couple $(\a_i,\b_i)$ we define the maps $F^i_1$, $F^i_2$ and $F^i_3$ as in the last section; we require that, if $\hat G_0$ is the solid triangle $ABC$, then points 1) and 2) at the end of section 1 hold; as we saw above, this is true if (2.5) holds. Now we can define the maps $F_\listo$ as in (1.22) and the set 
$\tilde G$ as in (1.24). 

As in (7), we set ${\cal R}=\{ (\a_i,\b_i) \}_{i\ge 1}$ and 
$\s{\cal R}=\{ (\a_{i+1},\b_{i+1}) \}_{i\ge 1}$. For $(\a_i,\b_i)$ as above, we set 
$$\e_i=\a_i\cdot\frac{5}{3},\qquad \l_i=\frac{3}{5}\e_i^2,\qquad
\tilde\l_{\rc,l}=\prod_{i=1}^l\l_i,\qquad
\tilde\e_{\rc,s}^l=\prod_{i=s}^l\e_i   .   \eqno (3.1)$$

Let $\Omega^s$ be as in (6) and let the pre-fractals $G_l$ be defined as after (4); for 
$\g_j$ defined as in (8) and $b>0$ to be determined later we define an energy on cables by  
$$\ec^2_{\{ \Omega^j \},s}(\nabla u,\nabla v)=
\frac{b}{1-\e(\Omega^s)}\sum_{j\in\Omega^s}
\int_0^1\diff u\circ\g_j(t)\cdot \diff v\circ\g_j(t)  \dr t   \eqno (3.2)$$ 
where $\e(\Omega^s)=\e_s$.

For a number $a>0$ to be determined later we define 
$$\hat\ec_{\rc,l}^1(\nabla u,\nabla v)=\frac{a}{\tilde\l_{\rc,l}}
\sum_{\listo\in\{ 1,2,3 \} }\sum_{j\in\{ AB,BC,AC \} }
\int_0^1\diff u\circ F_\listo\circ\g_j(t)\cdot
\diff v\circ F_\listo\circ\g_j(t)\dr t    \eqno (3.3)$$
and
$$\hat\ec_{\rc,l}^2(\nabla u,\nabla v)=
\sum_{s=1}^l
\frac{1}{\tilde\l_{\rc,s-1}}\cdot\frac{1}{\tilde\e^l_{\rc,s}}\cdot
\sum_{i_1,\dots,i_{s-1}\in\{ 1,2,3 \}}
\ec^2_{\{ \Omega^j \},s}(\nabla (u\circ F_{i_1,\dots,i_{s-1}}),\nabla (v\circ F_{i_1,\dots,i_{s-1}}))     
\eqno (3.4)$$
where $F_{i_1,\dots,i_{s-1}}$ is defined as in (1.22) and with the convention that, when $s=1$, $ F_{i_1,\dots,i_{s-1}}$ is the identity and $\tilde\l_{\rc,0}=1$. Clearly, in this way way we have the first equality below, while the second one comes from (3.2) since 
$\s\{ \Omega^j \}=\{ \Omega^{j+1} \}$ and $\e((\s\{ \Omega \})_s)=\e_{s+1}$. 
$$\hat\ec_{\s\rc,l}^2(\nabla u,\nabla v)=
\sum_{s=1}^l
\frac{1}{\tilde\l_{\s\rc,s-1}}\cdot\frac{1}{\tilde\e^l_{\s\rc,s}}\cdot
\sum_{i_1,\dots,i_{s-1}\in\{ 1,2,3 \}}
\ec^2_{\s\{ \Omega^j\} ,s }(\nabla (u\circ F_{i_1,\dots,i_{s-1}}),\nabla (v\circ F_{i_1,\dots,i_{s-1}}))=    $$
$$\sum_{s=1}^l
\frac{1}{\tilde\l_{\s\rc,s-1}}\cdot\frac{1}{\tilde\e^l_{\s\rc,s}}\cdot
\sum_{i_1,\dots,i_{s-1}\in\{ 1,2,3 \}}
\ec^2_{\{ \Omega^j \},s+1}(\nabla (u\circ F_{i_1,\dots,i_{s-1}}),\nabla (v\circ F_{i_1,\dots,i_{s-1}}))   .   $$

Lastly, we set 
$$\ec_{{\cal R},l}(\nabla u,\nabla v)=
\hat\ec_{\rc,l}^1(\nabla u,\nabla v)+\hat\ec_{\rc,l}^2(u,v)    .    
\eqno (3.5)$$
Just a word to justify the choice of the indices: when we are on the first pre-fractal and $l=1$, by (3.3) $\hat\ec^1_{\rc,1}$ contains the integrals on the three small triangles of figure 1, while $\hat\ec^2_{\rc,1}$ contains the integrals on the first generation cables which connect them, i. e. the elements of $\Omega^1$. Loosely speaking, (3.3) is the energy on the $l$-th pre-fractal $\tilde G_l$ of the disconnected gasket $\tilde G$ of (1.24), while (3.4) is the energy of all cables up to generation $l$, that is on all the cables tying together the triangles of $\tilde G_l$.

The following lemma shows that the forms $\ec_{\rc,l}$ are related by an iteration.

\lem{3.1} Let the form $\ec_{\rc,l}$ be defined by (3.5) and let $\l_i$, $\tilde\e^l_{r,s}$ be as in (3.1). Then, we have that 
$$\ec_{\rc,l+1}(\nabla u,\nabla v)=
\frac{1}{\l_1}\sum_{i=1}^3\ec_{\s\rc,l}(\nabla(u\circ F^1_i), \nabla(v\circ F^1_i))+$$
$$\frac{1}{\tilde\e^{l+1}_{\rc,1}}\ec^2_{\{ \Omega^j \},1}(\nabla u,\nabla v)  .    \eqno (3.6)$$

\proof We skip the proof of the following formula, which follows easily from (3.3) and is related to (1.27). 
$$\tilde\ec^1_{\rc,l+1}(\nabla u,\nabla v)=
\frac{1}{\l_1}
\sum_{i=1}^3\tilde\ec^1_{\s\rc,l}(\nabla (u\circ F^1_i),\nabla (v\circ F^1_i) )  .   $$

We show the recurrence for $\hat\ec^2_{\rc,l}$. The first equality below comes from (3.4); the second one comes from the formula after (3.4) and the fact that, by (3.1), 
$\l_1\cdot\tilde\l_{\s\rc,s}=\tilde\l_{\rc,s+1}$ and $\tilde\e^l_{\s\rc,s}=\tilde\e^{l+1}_{\rc,s+1}$. The third equality comes changing the indices while the last one is (3.4). 
$$\frac{1}{\l_1}\sum_{i=1}^3\hat\ec^2_{\s\rc,l}(\nabla (u\circ F^1_i),\nabla (v\circ F^1_i))=$$
$$\frac{1}{\l_1}\sum_{i=1}^3\sum_{s=1}^l
\frac{1}{\tilde\l_{\s\rc,s-1}}\cdot\frac{1}{\tilde\e^l_{\s\rc,s}}\cdot  
\sum_{i_1\dots i_{s-1}\in\{ 1,2,3 \}}
\ec^2_{\s\{ \Omega^j \},s}(\nabla(u\circ F^1_i\circ F^2_{i_1}\circ\dots\circ F^{s}_{i_{s-1}}),
\nabla(v\circ F^1_i\circ F^2_{i_1}\circ\dots\circ F^{s}_{i_{s-1}}))=$$
$$\sum_{s=1}^l
\frac{1}{\tilde\l_{\rc,s}}\cdot\frac{1}{\tilde\e^{l+1}_{\rc,s+1}}\cdot  
\sum_{i,i_1\dots i_{s-1}\in\{ 1,2,3 \}}
\ec^2_{\{ \Omega^j \},s+1}(\nabla(u\circ F^1_i\circ F^2_{i_1}\circ\dots\circ F^{s}_{i_{s-1}}),
\nabla(v\circ F^1_i\circ F^2_{i_1}\circ\dots\circ F^{s}_{i_{s-1}}))=$$
$$\sum_{s=2}^{l+1}
\frac{1}{\tilde\l_{\rc,s-1}}\cdot\frac{1}{\tilde\e^{l+1}_{\rc,s}}\cdot  
\sum_{i_1\dots i_{s-1}\in\{ 1,2,3 \}}
\ec^2_{\{ \Omega^j \},s}(\nabla(u\circ F^1_{i_1}\circ\dots\circ F^{s-1}_{i_{s-1}}),
\nabla(v\circ F^1_{i_1}\circ\dots\circ F^{s-1}_{i_{s-1}}))=$$
$$\hat\ec^2_{\rc,l+1}(\nabla u,\nabla v)-
\frac{1}{\tilde\e^{l+1}_{\rc,1}}\ec^2_{\{ \Omega^j \},1}(\nabla u,\nabla v)  .  $$

\fin

Now that we have an energy, we can we define harmonicity.

\vskip 1pc

\noindent{\bf Definition.} Let $G_l$ be defined as in (4) and let $\ec_{\rc,l}$ be as in (3.5). We say that $\ec_{\rc,l}$ is harmonic on $G_l$ if for all $u\in C^2(\R^2,\R)$ there is a bounded Borel function $\fun{g}{G_l}{\R}$ such that, for all $v\in C^1(\R^2,\R)$ with $v(A)=v(B)=v(C)=0$ we have 
$$\ec_{{\cal R},l}(\nabla u,\nabla v)=
-\int_{G_l}g(x)v(x)\dr H^1(x)  \eqno (3.7)$$
where $H^1$ denotes the one-dimensional Haussdorff measure on $G_l$. 

We say that $u\in C^1(\R^2,\R)$  is harmonic on $G_l$ if, for every test function $v$ as above, we have 
$$\ec_{{\cal R},l}(\nabla u,\nabla v)=0  \eqno (3.8)$$

\vskip 1pc

\noindent{\bf Remarks.} By (3.5), $\ec_{{\cal R},l}$ is a sum of integrals on the edges of the graph $G_l$; formula (3.7) says that, if we integrate by parts on the edges of $G_l$, the only boundary terms which survive are those in $A$, $B$, $C$ and they too disappear if 
$v$ vanishes on them. In a sense, we are stipulating that  the curvature terms of the Laplace-Beltami operator, which concentrate on the vertices since the edges are straight lines, are zero. 

The name "harmonic" for the pre-fractals satisfying (3.7) comes from the analogous property of the pre-fractals of the harmonic Sierpinski gasket, the first of which is depicted in figure 2 above. 

We must resist the temptation to call the function $g$ of (3.7) the "Laplacian" of $u$; in fact, our "Laplacian" is a different function because the measure with which we endow 
$G_l$ in (6.14) below, though absolutely continuous with respect to $H^1$, does not coincide with $H^1$: the "Laplacian" will be the function $g$ multiplied by some density. Naturally, this has no effect when $u$ is harmonic, i. e. when $g\equiv 0$. 

Keeping in mind the last remark, the integration by parts formula on $G_l$ will eventually imply an integration by parts formula on the full fractal $G$, as in [19]. 

\lem{3.2} Let us suppose that the form $\ec_{{\cal R},l}$ of (3.5) is harmonic on $G_l$, and let $\fun{u}{\R^2}{\R}$ be affine. Then, $u$ is harmonic on $G_l$, i. e. formula (3.8) holds. Conversely, if all affine functions are harmonic, then (3.7) holds. 

\proof We skip the proof of the converse, since it is the same argument we use in section 7 below to prove the integration by parts formula. 

The argument for the direct part is very simple: if we integrate by parts in (3.5) to get the Laplacian, by (3.7) we find no terms which concentrate on the vertices of $G_l$; as for the edges, no terms concentrate on them because the functions 
$u\circ F_{i_1\dots i_l}\circ\g_s$ are affine and thus their second derivative vanishes. Now to the rigorous proof. 

Let $\tilde G$ be the disconnected gasket of (1.24) and let $\tilde G_l$ be the $l$-th pre-fractal of $\tilde G$, i. e. 
$$\tilde G_l=
\bigcup_{{i_1\dots i_l}\in \{ 1,2,3 \}}F_{i_1\dots i_l}(G_0)$$
where $G_0$ is the triangle of (4) and $F_{i_1\dots i_l}$ is defined as in (1.22). 

Let us look at the affine function $u$ in the integrals of (3.5). The form $\ec_{{\cal R},l}$ is a sum of two terms: the first one is $\hat\ec_{\rc,l}^1$ and it contains by (3.3) the integrals on the edges of $\tilde G_l$; using the notation of (8) and (1.22), $u$ appears in them as 
$$\fun{u\circ F_\listo\circ\g_j}{[0,1]}{\R}$$
with 
$$i_1,\dots, i_l\in\{ 1,2,3 \},\qquad j\in \{ AB,BC,AC \}   .   $$
The second term is (3.4) and it contains the integrals on the "cables"; with the usual convention that $F_{i_1\dots i_0}=id$, $u$ appears in them as 
$$\fun{u\circ F_{i_1\dots i_r}\circ\g_z}{[0,1]}{\R}$$
with 
$${i_1\dots i_r}\in\{ 1,2,3 \},\qquad 
1\le r\le l-1,\qquad
z\in \Omega^{r+1}   .   $$
In both cases, these functions of $t\in[0,1]$ are affine, since they are compositions of affine functions. In particular , 
$$\frac{\dr^2}{\dr t^2}u\circ F_{i_1\dots i_l}\circ\g_s\equiv
\frac{\dr^2}{\dr t^2}u\circ F_{i_1\dots i_r}\circ\g_z\equiv 0  .  $$
Thus, if we integrate by parts in (3.5), the only terms which survive are the boundary ones, which are calculated at the vertices of $G_l$: if we show that the sum of the boundary terms at all vertices is zero, the lemma follows. 

Let $S$ be a vertex of $G_l$, i. e. $S=F_{i_1\dots i_l}(P)$, with $P\in\{ A,B,C \}$. Inspection of figure 1 shows the following: if $S\not\in\{ A,B,C \}$, then at $S$ lands exactly one cable, $F_{i_1\dots i_{r-1}}\circ\g_z(h)$ with $h\in\{ 0,1 \}$, $z\in\Omega^r$ and 
$1\le r\le l$. Moreover, at $S$ land two edges of the pre-fractal $\tilde G_l$, namely the two segments $F_\listo\circ\g_s$ such that $F_\listo\circ\g_s(b)=S$ for some $b\in\{ 0,1 \}$ and $s\in\{ AB,BC,AC \}$. By formulas (3.3), (3.4) and (3.5) this implies that the sum of the boundary terms at $S$ is given by 
$$-v(S)\cdot\left[
\frac{a}{\tilde\l_{\rc,l}}\sum_{(s,b)}\diff u\circ F_{i_1\dots i_l}\circ\g_s(b)(-1)^b+
\frac{ b}{\tilde\l_{\rc,r-1}\tilde\e^l_{\rc,r}(1-\e_r)}
\diff u\circ F_{i_1\dots i_{r-1}}\circ\g_{z}(h)(-1)^h
\right]     \eqno (3.9)$$
where the sum is on the two edges of $\tilde G_l$ we mentioned above. Formula (3.7) tells us that all the terms in square parentheses vanish, save possibly those at $A$, $B$, $C$. These too, however, are zero. Indeed, by figure 1 no cable arrives at $A$, $B$ or $C$; thus, for $S\in\{ A,B,C \}$ we have a similar formula to (3.9) but without the cables: for instance, in $A$, 
$$-v(A)\cdot\frac{1}{\tilde\l_{\rc,l}}
\sum_s\diff u\circ F^1_1\circ\dots\circ F^l_1\circ\g_s(0)  .   \eqno (3.10)$$
The sum is on the two elements $s\in \{ AB,AC \}$ which satisfy 
$F^1_1\circ\dots\circ F^l_1\circ\g_s(0)=A$. Now the expression above is zero because we are stipulating that the test function $v$ satisfies $v(A)=0$. 

\fin

\vskip 2pc

\centerline{\bf \S 4}
\centerline{\bf Harmonicity on $G_1$}

\vskip 1pc

Let $G_l$ be the set defined in (4). In particular, 
$$G_1=\bigcup_{i=1}^3 F_i(G_0)\cup\Omega^1   \eqno (4.1)$$
where $\Omega^1$ is as (6) and $G_0$ is as in (4); the graph $G_1$ is homeomorphic to that of figure 1. 

We are going to find conditions on the maps $F_1^1$, $F_2^1$, $F_3^1$ such that 
$\ec_{{\cal R},1}$ satisfies (3.7) on $G_1$. In other words, if we integrate by parts on the edges of $G_1$, we want the boundary terms to cancel out, save those in $A$, $B$, $C$. Inspection of (4.1) (or of figure 1, since $G_1$ has the same topology) shows that there are six vertices of $G_1$ different from $A$, $B$, $C$: 
$$F_1^1(B), \quad F_1^1(C),\quad F_2^1(A),\quad F_2^1(C),\quad F_3^1(A),\quad F_3^1(B)  .  $$
By symmetry, it suffices to impose that the sum of the boundary terms at $F_1^1(B)$ vanishes. We list the edges that land at this point; with the notation of (8) and looking again at figure 1, they are 
$$F_1^1\circ\g_{BA},\qquad F_1^1\circ\g_{BC},\qquad \g_{F_1^1(B)F_2^1(A)}  .  $$
Note that, by (8), the three curves above are in $F_1^1(B)$ at $t=0$; as a consequence, we have that, in (3.9), $(-1)^b=(-1)^h=1$; moreover, $l=r=1$. 

Let $\e_1$ be as in (3.1). On the graph $G_1$ there is the form $\ec_{{\cal R},1}$ of (3.5); we integrate by parts in this formula. The only contributions to the boundary term at $F_1(B)$ are those of the three edges above; we apply the chain rule to (3.9) and recall that $\tilde\l_{\rc,0}=1$; we get that (3.9) vanishes at $S=B$ for all $u,v\in C^1(\R^2,\R)$ if and only if 
$$\frac{a}{\l_1}\left[
\diff F_1^1\circ\g_{BA}(0)+\diff F_1^1\circ\g_{BC}(0)
\right]+\frac{b}{\e_1(1-\e_1)}\diff\g_{F_1^1(B)F_2^1(A)}(0)=0  .   \eqno (4.2)$$
Since 
$$DF_1^1=T_1=\left(
\matrix{
\a_1,&0\cr
0,&\b_1
}
\right)  ,  $$
the chain rule and (2.1), (2.8), (2.9), (2.10) imply that (4.2) is equivalent to the equality below. 
$$0=\frac{a}{\l_1}
\left(
\matrix{
\a_1,&0\cr
0 &\b_1
}
\right)
\left[
\left(
\matrix{
-\frac{\sqrt 3}{2}\cr
-\2
}
\right) +
\left(
\matrix{
0\cr
-1
}
\right)
\right]  +
\frac{b}{\e_1(1-\e_1)}
\left(
\matrix{
(2-\b_1-3\a_1)\cdot\frac{\sqrt 3}{4}\cr
(2-\b_1-3\a_1)\cdot\frac{1}{4}
}
\right)      .  $$
The last equation is equivalent to the system 
$$\left\{
\eqalign{
\frac{b}{\e_1(1-\e_1)}(2-\b_1-3\a_1)&=\frac{2\a_1 a}{\l_1}\cr
\frac{b}{\e_1(1-\e_1)}(2-\b_1-3\a_1)&=\frac{6\b_1 a}{\l_1}  .  
}
\right.  $$
Note that, by the second formula of (2.5), $2-\b_1-3\a_1>0$;  thus, for al $a>0$ there is $b>0$ which satisfies the formula above if and only if  
$$\a_1=3\b_1  .   \eqno (4.3)$$

This brings us to the following lemma.  

\lem{4.1} Let 
$$(\a_i,\b_i)=\e_i\left(
\frac{3}{5},\frac{1}{5}
\right),\qquad   \e_i\in (0,1)    \eqno (4.4)$$
and let the maps $F^i_1$, $F^i_2$, $F^i_3$ be defined as in section 2 for $(\a_i,\b_i)$ as above. Then, the following holds. 

\noindent 1) Conditions 1) and 2) before (1.22) hold for the solid triangle $\hat G_0$, which we defined after (2.4). 

\noindent 2) If $b=a$, then the pre-fractal $G_1$ is harmonic. 

\noindent 3) For all $i\ge 1$ the triple of maps $(F_1^i, F_2^i, F_3^i)$ defined in section 2 satisfies hypothesis (ND) of section 1. 

\proof Point 2) has been proven at the beginning of this section: by the system before (4.3) we get that 
$$b=\frac{3}{5}\cdot
\frac{\e_1^2a}{\l_1}=
a        $$
where the last equality comes if we take $\l_1$ as in (3.1).   

As for point 1), we call $F_1$, $F_2$, $F_3$ the maps of (2.1), (2.3) and (2.4) for 
$(\a,\b)=\left( \frac{3}{5},\frac{1}{5} \right)$, i. e. for $\e_i=1$; these maps are well-studied, since they yield the harmonic gasket of figure 2. We continue to call $\hat G_0$ the solid triangle $ABC$ and we recall from [9] that $F_1(\hat G_0)$, $F_2(\hat G_0)$ and $F_3(\hat G_0)$ intersect only at the three points $a,b,c$ of figure 2. Thus, conditions 1) and 2) before (1.22) follow if we show the following: when $\e_i\in(0,1)$, 
$F^i_1(\hat G_0)$ is contained in $F_1(\hat G_0)\setminus\{ b,c \}$, 
$F^i_2(\hat G_0)$ is contained in $F_2(\hat G_0)\setminus\{ a,c \}$ and 
$F^i_3(\hat G_0)$ is contained in $F_3(\hat G_0)\setminus\{ a,b \}$. By symmetry, it suffices to verify the assertion for $F^i_1$, but this follows easily from (2.1).  

Point 3) follows from the three points below. 

\noindent $\bullet$) As shown in [4], (ND) holds for the harmonic Sierpinski gasket, where $\a=\frac{3}{5}$ and $\b=\frac{1}{5}$.

\noindent $\bullet$) Let $\e>0$; it is immediate that (ND) holds for the derivatives $DF_1$, $DF_2$, $DF_3$ if and only if it holds for $\e DF_1^1$, $\e DF_2^1$, $\e DF_3^1$. 

\noindent $\bullet$) By formula (4.4),  
$$\frac{1}{\e_i}DF_1^i=\left(
\matrix{
\frac{3}{5},&0\cr
0,&\frac{1}{5}
}
\right)    ,    $$
i. e. the derivative of the map $F^i_1$ (and consequently that of $F^i_2$ and $F^i_3$) is a multiple of the derivative of the corresponding map for the harmonic gasket. 

\fin

\vskip 2pc

\centerline{\bf \S 5}

\centerline{\bf The $l$-th pre-fractal is harmonic}

\vskip 1pc

In lemma 5.1 below we are going to show that, if all couples $(\a_i,\b_i)$ satisfy (4.4), then $G_l$ is harmonic; before proving it, we need to fix some notation. We suppose that 
$(\a_i,\b_i)$ are defined as in (4.4); as in (7) we set 
$\rc=\left\{ \e_i\left( \frac{3}{5},\frac{1}{5}\right) \right\}_{i\ge 1}$ and 
$\s\rc=\left\{ \e_{i+1}\left( \frac{3}{5},\frac{1}{5}\right) \right\}_{i\ge 1}$. We re-write (3.1) under this condition. 
$$\l_i=\frac{3}{5}\e_i^2,\qquad
\tilde\l_{\rc,l}=\prod_{i=1}^l\l_i,\qquad
\tilde\e^l_{\rc,s}=\prod_{i=s}^l\e_i  .  \eqno (5.1)$$

We explain the choice of the constants $\l_i$ in (5.1); since $\l_i$ is the eigenvalue of a Ruelle operator, we must expand a little on this topic. 

We define a matrix-valued measure $\tau_0$ supported on the set $G_0$ of (4) in the following way: let $E\in C(\bar B(0,R),M)$ be an arbitrary, continuous field of symmetric matrices, let $P_v$ denote the orthogonal projection on the vector $v\in\R^2$; for  $a>0$ to be determined presently we set 
$$\int_{\bar B(0,R)}(E,\dr\tau_0)_{HS}=$$
$$a\sum_{j\in\{ AB,BC,AC \}}
\int_0^1(E_{\g_j(t)},P_{\dot\g_j(t)})_{HS}||\dot\g_j(t)||^2\dr t  .  $$
Now note that the eigenvector $Q$ of theorem 1.1 is positive-definite (by lemma 4.1 this is the identity, as for the harmonic gasket); thus, if in the formula above we set $E=Q$, each of the three integrals in the sum on the right is positive. In particular, we can find $a>0$ such that $\tau_0\in\pc_Q(\bar B(0,R),M)$ (a space we defined before (1.18)) and this will be our choice. 

Note that, with the definition above, if $u,v\in C^1(\R^2,\R)$, then the chain rule and second equality of (1.12) imply that 
$$\int_{\bar B(0,R)}(\nabla u,\dr\tau_0\cdot\nabla v)=$$
$$a\sum_{j\in\{ AB,BC,AC \}}
\int_0^1\diff u\circ\g_j(t)\cdot\diff v\circ\g_j(t)\dr t  .  $$

For the maps $F^i_1$, $F^i_2$, $F^i_3$ determined by the constants $(\a_i,\b_i)$ of (4.4) we define the Ruelle operators $\L_i$ on the coding space as after (1.24). 
$$\fun{\L_{i}}{C(\Sigma,M)}{C(\Sigma,M)}$$
$$\L_{i}(A)=\sum_{j=1}^3 \tr(DF^i_j)\cdot A_{(jx)}\cdot (DF^i_j)   .   $$
The operator $\fun{\L_\Sigma}{C(\Sigma,M)}{C(\Sigma,M)}$ is defined as above, but with 
$\e_i=1$; it is the Ruelle operator on the Harmonic Sierpinski Gasket  (or better, on its coding) and it is well-known ([4], [9]) that its eigenvalue $\l$ is equal to $\frac{3}{5}$. By the arguments before (1.25), this implies that, for all $i\ge 1$, 
$$\frac{1}{\l_i}\L_{i}=\frac{1}{\l}\L_\Sigma   
\txt{and}  \l_i=\e_i^2\l=
\frac{3}{5}\e_i^2  ,  \eqno (5.2)$$
which explains the first equality of (5.1). Analogously, we define $\tilde\L_{\tilde G,l}$ as in (1.26) and $\tau_l$ as in (1.28). The definition of the adjoint and a comparison with (3.3) shows that 
$$\hat\ec^1_{\rc,l}(\nabla u,\nabla v)=
\int_{\bar B(0,R)}(\nabla u,\dr\tau_l\cdot\nabla v)   .    \eqno (5.3)$$

\lem{5.1} Let the coefficients $(\a_i,\b_i)$ satisfy (4.4) for a sequence 
$\{ \e_i \}_{i\ge 1}\subset(0,1)$; let $a>0$ be as at the beginning of this section; let $\l_i$ be the eigenvalue of $\L_{i}$ as in (5.2). Let $b=a$ as in lemma 4.1, and let $\tilde\l_{\rc,l}$, $\tilde\e^l_{\rc,s}$ be as in (5.1). Then, for all $l\ge 1$, the form 
$\ec_{\rc,l}$ of (3.5) is harmonic, i. e. it satisfies (3.7). 

\proof The proof is by induction, with lemma 4.1 catering for $G_1$. 

Let us suppose that $G_l$ is harmonic; we prove that $G_{l+1}$ is harmonic. 

We re-write the recurrence (3.6), using (3.2). 
$$\ec_{{\cal R},l+1}(\nabla u,\nabla v)=
\frac{1}{\lambda_1}\sum_{j=1}^3
\ec_{\s{\cal R},l}(\nabla u\circ F_j^1,\nabla v\circ F_j^1)+    \eqno (5.4)_a$$
$$\frac{ b}{\tilde\e^{l+1}_{\rc,1}(1-\e_1)}\cdot\sum_{j\in\Omega^1}
\int_0^1\diff u\circ\g_j(t)\cdot \diff u\circ\g_j(t)\dr t  .  \eqno (5.4)_b$$

First of all, by the induction hypothesis, the boundary terms of 
$\ec_{\s\rc,l}(\nabla u\circ F_i^1,\nabla v\circ F_i^1)$ are all zero, save the ones at $A$, 
$B$, $C$; by $(5.4)_a$ they contribute to the boundary terms of 
$\ec_{\rc,l+1}(\nabla u,\nabla v)$ at the following nine vertices of $G_{l+1}$: 
$$F_i^1(A),\qquad F_i^1(B),\qquad F_i^1(C),\qquad i\in\{ 1,2,3 \}  .  $$
We integrate by parts in $(5.4)_b$; by (6), the boundary terms concentrate on a subset of the set above, i. e. on the following six vertices: 
$$F_1^1(B),\qquad F_1^1(C),\qquad F_2^1(A),\qquad F_2^1(C),\qquad 
F_3^1(A),\qquad F_3^1(B) .   \eqno (5.5)$$

Let us show that the boundary terms vanish at the nine vertices above. The three points $A=F_1^1(A)$, $B=F_2^1(B)$, $C=F_3^1(C)$ are dispatched easily: the boundary terms at these points vanish by (3.10) because we are stipulating that $v(A)=v(B)=v(C)=0$. 

We show that the boundary terms vanish at the remaining six points, those of (5.5). By symmetry, it suffices to show that they vanish at $F_1^1(B)$. 

We integrate by parts in the left hand side of $(5.4)_a$; by (3.5), or directly by (3.9), we see that the boundary term of $\ec_{{\cal R},l+1}(\nabla u,\nabla v)$ at $F^1_1(B)$ is, up to the sign, 
$$v(F^1_1(B))\cdot\Big[
\frac{a}{\tilde\l_{\rc,l+1}}\diff u\circ F^1_1\circ F^2_2\circ\dots\circ F^{l+1}_2\circ\g_{BA}(0)+$$
$$\frac{a}{\tilde\l_{\rc,l+1}}\diff u\circ F^1_1\circ F^2_2\circ\dots\circ F^{l+1}_2\circ\g_{BC}(0)+
\left(    \frac{a}{\tilde\e^{l+1}_{\rc,1}(1-\e_1)}   \right) 
\diff u\circ\g_{F^1_1(B)F^1_2(A)}(0)
\Big]   .  $$
Just a word of explanation: the first two summands come from integration by parts in (3.3), namely from the two edges of $\tilde G_{l+1}$ which contain $F^1_1(B)$; since  $F^j_2(B)=B$ for all $j$, these edges are 
$F^1_1\circ F^2_2\circ\dots\circ F^{l+1}_2\circ\g_{BA}$ and 
$F^1_1\circ F^2_2\circ\dots\circ F^{l+1}_2\circ\g_{BC}$. The last summand comes from $(5.4)_b$, namely from the element of $\Omega^1$ which contains $F^1_1(B)$. 

We have to show that the expression above is zero for all $v,u\in C^1(\R^2,\R)$; by the chain rule this is tantamount to the following vector equality. 
$$\frac{a}{\tilde\l_{\rc,l+1}}\diff  F^1_1\circ F^2_2\circ\dots\circ F^{l+1}_2\circ\g_{BA}(0)+$$
$$\frac{a}{\tilde\l_{\rc,l+1}}\diff F^1_1\circ F^2_2\circ\dots\circ F^{l+1}_2\circ\g_{BC}(0)+
\left(    \frac{a}{\tilde\e^{l+1}_{\rc,1}(1-\e_1)}   \right) 
\diff \g_{F^1_1(B)F^1_2(A)}(0)  =0  .  \eqno (5.6)$$
Formula (4.2) tells us that 
$$\frac{a}{\l_1}\diff  F^1_1\circ\g_{BA}(0)+
\frac{a}{\l_1}\diff  F^1_1\circ\g_{BC}(0)+
\left(    \frac{a}{\e_1(1-\e_1)}   \right) 
\diff \g_{F^1_1(B)F^1_2(A)}(0)  =  0   .  $$
Let us compare the last two formulas: if we multiply the second one by 
$\frac{\e_1}{\tilde\e^{l+1}_{\rc,1}}$ we get that (5.6) holds if 
$$\frac{1}{\tilde\l_{\rc,l+1}}\diff  F^1_1\circ F^2_2\circ\dots\circ F^{l+1}_2\circ\g_{BA}(0)+
\frac{1}{\tilde\l_{\rc,l+1}}\diff  F^1_1\circ F^2_2\circ\dots\circ F^{l+1}_2\circ\g_{BC}(0) =$$
$$\frac{\e_1}{\tilde\e^{l+1}_{\rc,1}}\cdot
\left[
\frac{1}{\l_1}\diff  F^1_1\circ\g_{BA}(0)+
\frac{1}{\l_1}\diff  F^1_1\circ\g_{BC}(0)
\right]   ,    $$
where we have simplified the constant $a>0$. By the chain rule and recalling that $DF^1_1$ is invertible, the last formula holds if  
$$\frac{1}{\tilde\l_{\rc,l+1}}D(F_2^2\circ\dots\circ F^{l+1}_2)[\dot\g_{BA}(0)+\dot\g_{BC}(0)]=
\frac{1}{\l_1}\cdot\frac{\e_1}{\tilde\e^{l+1}_{\rc,1}}
[\dot\g_{BA}(0)+\dot\g_{BC}(0)]  .  \eqno (5.7)$$
Formulas (1) and (8) (or, if you prefer, direct inspection of figure 1) show that, for some 
$\d>0$, 
$$\dot\g_{BA}(0)+\dot\g_{BC}(0)=w$$
with 
$$w=-\d\left(
\matrix{
\2\cr
\frac{\sqrt 3}{2}
}
\right)  .   $$
We re-write (5.7) with this notation, getting the equivalent formula 
$$\frac{1}{\tilde\l_{\rc,l+1}}\cdot DF^2_2\cdot\dots\cdot DF^{l+1}_2 w=
\frac{1}{\l_1}\cdot \frac{\e_1}{\tilde\e^{l+1}_{\rc,1}}\cdot w   .  \eqno (5.8)$$
An easy calculation using (4.4) and the formula for $T^i_2=DF^i_2$ before (2.6) shows that 
$$DF^i_2 w=\e_i\cdot\frac{3}{5}\cdot w  .  $$
This implies the second equality below; as for the first one, it comes from (5.2), which implies that $\l_i=\frac{3}{5}\e_i^2$, and the definition of $\tilde\e^{l+1}_{\rc,1}$ in (5.1). 
$$\frac{\tilde\e^{l+1}_{\rc,1}}{\tilde\l_{\rc,l+1}}\cdot
DF^2_2\cdot\dots\cdot DF^{l+1}_2 w=$$
$$\frac{1}{\prod_{i=1}^{l+1}\e_i}\cdot
\frac{1}{\left( \frac{3}{5} \right)^{l+1}}DF^2_2\cdot\dots\cdot DF^{l+1}_2 w=
\frac{1}{\frac{3}{5}\e_1}w  .  \eqno (5.9) $$
Since 
$$\frac{\e_1}{\l_1}w=\frac{1}{\frac{3}{5}\e_1}w$$
by (5.2), formula (5.8) follows.

\fin

\vskip 2pc

\centerline{\bf \S 6}

\centerline{\bf The limit form and the natural measure}

\vskip 1pc

Let $\hat\ec^1_{\rc,l}(\nabla u,\nabla v)$ be as in (3.3) for a sequence 
$(\a_i,\b_i)$ which satisfies (4.4); the numbers $\l_i$ and $a$ are as in the last section. 
We recall formula (1.29): for all $u,v\in C^1(\R^2,\R)$, 
$$\hat\ec^1_{\rc,l}(\nabla u,\nabla v)\tends\ec^1_\rc(\nabla u,\nabla v)
\txt{as}l\tends+\infty    \eqno (6.1)$$
where 
$$\ec^1_\rc(\nabla u,\nabla v)=
\int_{\bar B(0,R)}(\nabla u,\dr\tau_\rc\cdot\nabla v)  .  $$

Let $\tau_\Sigma$ and $\kappa_\Sigma$ be as at the end of section 1; by point 3) of theorem 1.1 these two measures are mutually absolutely continuous with bounded Radon-Nikodym derivatives; since $\Phi$ is bijective, the push-forwards $\tau$ and $\kappa$ share the same property and are supported on $G$. In particular, we can write 
$\tau=T_x\kappa$ where 
$\fun{T}{G}{M}$ is a bounded Borel function. With this notation, we get 
$$\ec^1_\rc(\nabla u,\nabla v)=
\int_G(T_x\nabla u(x),\nabla v(x))\dr\kappa(x)    .     \eqno (6.2)$$

Let now $\ec^2_{\rc,l}$ be as in (3.4): it is the part of the energy carried by the cables up to generation $l$. We want to show that $\ec^2_{\rc,l}$ too converges. 

\vskip 1pc

\lem{6.1} Let $(\a_i,\b_i)$ be as in (4.4), let $\tilde\l_{\rc,r},\tilde\e^l_{\rc,s}$ be as in (5.1) and let $\Omega^s$ be as in (6); let $b>0$ be as in lemma 4.1. We suppose that 
$$\prod_{i\ge 1}\e_i>0  .  \eqno (6.3)$$
In particular, for some $\d>0$ independent of $s$ we have the first inequality below; the equality is the definition of $\tilde\e^\infty_{\rc,s}$ and the second inequality comes from the fact that $\e_i<1$ for all $i$.  
 $$\d\le\tilde\e^\infty_{\rc,s}=\prod_{i=s}^{+\infty}\e_i < 1    
 \txt{for all $s\ge 1$. }   \eqno (6.4)$$
For $u,v\in C^1(\R^2,\R)$ we consider the sum
$$\ec^2_\rc(\nabla u,\nabla v)=$$
$$\sum_{s\ge 1}\sum_{j\in\Omega^s}\sum_{i_1\dots i_{s-1}\in \{ 1,2,3 \} }
\frac{a}{\tilde\l_{\rc,s-1}\tilde\e^\infty_{\rc,s}(1-\e_s)}\cdot
\int_0^1\diff u\circ F_{i_1\dots i_{s-1}}\circ\g_j(t)\cdot 
\diff v\circ F_{i_1\dots i_{s-1}}\circ\g_j(t)  \dr t   .    \eqno (6.5)$$
We continue in the convention that, when $s=0$, $F_{i_1\dots i_s}$ is the identity and 
$\tilde\l_{\rc,0}=1$. We assert that the sum on the right converges absolutely for all 
$u,v\in C^1(\R^2,\R)$. In particular,  $\ec^2_\rc$ is defined on $C^1(\R^2,\R)$.  

Moreover, for all $u,v\in C^1(\R^2,\R)$ and the form $\ec^2_{\rc,l}$ of (3.4) we have that 
$$\ec_{\rc,l}^2(\nabla u,\nabla v)\tends\ec^2_\rc(\nabla u,\nabla v)
\txt{as}l\tends+\infty  .  \eqno (6.6)$$

\proof The second inequality below follows from the fact that, by (5.1), the product on the right contains more factors $\e_i\in(0,1)$ than the product on the left.  
$$0<\frac{1}{\tilde\e^l_{\rc,s}}\le\frac{1}{\tilde\e^\infty_{\rc,1}}\qquad
\forall l,s\ge 1  .  $$
In particular, for fixed $s$ the sequence $\{ \frac{1}{\tilde\e^l_{\rc,s}} \}_{l\ge 1}$ is  bounded, increasing in $l$ and converging to $(\tilde\e^\infty_{\rc,s})^{-1}$ as 
$l\tends+\infty$; these numbers are uniformly bounded in $s$ by (6.4). Thus, each term of the sum (3.4) converges to the corresponding term of (6.4); by dominated convergence for  series, (6.5) follows if we show that  
$$\sum_{s\ge 1}\sum_{j\in\Omega^s}\sum_{i_1\dots i_{s-1}\in \{ 1,2,3 \} }
\frac{ 1}{\tilde\l_{\rc,s-1}(1-\e_s)}\cdot\left\vert
\int_0^1\diff u\circ F_{i_1\dots i_{s-1}}\circ\g_j(t)\cdot 
\diff v\circ F_{i_1\dots i_{s-1}}\circ\g_j(t)  \dr t   
\right\vert    <+\infty  .    $$

Now recall that $u$ and $v$ have bounded derivatives in a ball $B(0,R)$ which contains all the sets $G_l$. Using this fact and the chain rule we get that the formula above is implied by 
$$\sum_{s\ge 1}\sum_{j\in\Omega^s}\sum_{i_1\dots i_{s}\in \{ 1,2,3 \} }
\frac{1}{\tilde\l_{\rc,s-1}(1-\e_s)}\cdot
\int_0^1\left\vert
\diff F_{i_1\dots i_{s-1}}\circ\g_j(t)
\right\vert^2  \dr t
<+\infty  .  \eqno (6.7)$$
Before showing (6.7) we recall from [4] that $\tau(\tilde G)$ and $Q$ are positive-definite matrices (actually, multiples of the identity). Now the Gibbs formula, i. e. point 6) of theorem 1.1, implies that, for some $C_1>0$ independent of $i_1\dots i_{s}$ and $s$, 
$$\frac{1}{\tilde\l_{\rc,s-1}}
||DF_{i_1\dots i_{s-1}}||_{HS}^2\le C_1\kappa([i_1\dots i_{s-1}]_{\tilde G})
\eqno (6.8)$$
where the notation is that of (1.22). Moreover, from (6) and the definition of the triple $F^s_1,F^s_2,F^s_3$ we easily get that, for some $C_4>0$ independent of $s$, 
$$\left\vert
\diff\g_j(t)
\right\vert  \le C_4(1-\e_s)
\txt{if}j\in\Omega^s  .  \eqno (6.9)$$

We show (6.7): the first inequality below comes from the properties of the Hilbert-Schmidt norm, (6.9) and the fact that $\Omega^s$ has three elements; the second inequality comes from (6.8); the last equality comes from the fact that $\kappa$ is probability and the sets 
$[i_1\dots i_{s-1}]_{\tilde G}$ are disjoint, a fact we saw in section 1. The last inequality follows from (6.3), taking logarithms. 
$$\sum_{s\ge 1}\sum_{j\in\Omega^s}\sum_{i_1\dots i_{s-1}\in \{ 1,2,3 \} }
\frac{1}{\tilde\l_{\rc,s-1}(1-\e_s)}\cdot
\int_0^1\left\vert
\diff F_{i_1\dots i_{s-1}}\circ\g_j(t)
\right\vert^2  \dr t\le$$
$$C_3\sum_{s\ge 1}\sum_{i_1\dots i_{s-1}\in \{ 1,2,3 \} }
\frac{1}{\tilde\l_{\rc,s-1}(1-\e_s)}||DF_{i_1\dots i_{s-1}}||^2_{HS}(1-\e_s)^2\le$$
$$C_3C_1\sum_{s\ge 1}\sum_{i_1\dots i_{s-1}\in \{ 1,2,3 \} }
\kappa([i_1\dots i_{s-1}]_G)(1-\e_s)=
C_3C_1 \sum_{s\ge 1}(1-\e_s)<+\infty .  $$

\fin

For $\ec^1_\rc$ as in (6.2) and $\ec^2_\rc$ as in (6.5) we define 
$$\fun{\ec_\rc}{C^1(\R^2,\R)\times C^1(\R^2,\R)}{\R}$$
$$\ec_\rc(\nabla u,\nabla v)=\ec^1_\rc(\nabla u,\nabla v)+\ec^2_\rc(\nabla u,\nabla v)  .   
\eqno (6.10)$$
Now (6.1) and (6.6) imply that, for all $u,v\in C^1(\R^2,\R)$, 
$$\ec_{\rc,l}(\nabla u,\nabla v)\tends\ec_\rc(\nabla u,\nabla v)  .  \eqno (6.11)$$
We end this section re-writing the form $\ec_\rc$ of (6.10) in a way more similar to (6.2). 

We begin with (6.5), which we write as our first equality below; the second equality comes from the chain rule; for the third one we define $P_v$ as the orthogonal projection on the space generated by $v$. The fourth one comes from the definition of the push-forward of a measure; we have set 
$\tilde P_x=P_{DF_{i_1\dots i_{s-1}}\dot\g_j(t)}$ if $x=F_{i_1\dots i_{s-1}}\circ\g_j(t)$. The last equality is the definition of the measure $m$. 
$$\ec^2_\rc(\nabla u,\nabla v)=$$
$$\sum_{s\ge 1}\sum_{j\in\Omega^s}\sum_{i_1\dots i_{s-1}\in\{ 1,2,3 \}}
\frac{ a}{\tilde\l_{\rc,s-1}\tilde\e^\infty_{\rc,s}(1-\e_s)}\cdot
\int_0^1\diff u\circ F_{i_1\dots i_{s-1}}\g_j(t)\cdot
\diff v\circ F_{i_1\dots i_{s-1}}\g_j(t)\dr t  =  $$
$$\sum_{s\ge 1}\sum_{j\in\Omega^s}\sum_{i_1\dots i_{s-1}\in\{ 1,2,3 \}}
\frac{ a}{\tilde\l_{\rc,s-1}\tilde\e^\infty_{\rc,s}(1-\e_s)}\cdot$$
$$\int_0^1(\nabla u\vert_{F_{i_1\dots i_{s-1}}\circ\g_j(t)},DF_{i_1\dots i_{s-1}}\cdot\dot\g_j(t))\cdot 
(\nabla v\vert_{F_{i_1\dots i_{s-1}}\circ\g_j(t)},DF_{i_1\dots i_{s-1}}\cdot\dot\g_j(t))\dr t=$$
$$\sum_{s\ge 1}\sum_{j\in\Omega^s}\sum_{i_1\dots i_{s-1}\in\{ 1,2,3 \}}
\frac{ a}{\tilde\l_{\rc,s-1}\tilde\e^\infty_{\rc,s}(1-\e_s)}\cdot$$
$$\int_0^1||P_{DF_{i_1\dots i_{s-1}}\dot\g_j(t)}\nabla u|_{F_{i_1\dots i_{s-1}}\circ\g_j(t)}||\cdot
||P_{DF_{i_1\dots i_{s-1}}\dot\g_j(t)}\nabla v|_{F_{i_1\dots i_{s-1}}\circ\g_j(t)}||\cdot
||DF_{i_1\dots i_{s-1}}\dot\g_j(t)||^2\dr t  =  $$
$$\sum_{s\ge 1}\sum_{j\in\Omega^s}\sum_{i_1\dots i_{s-1}\in\{ 1,2,3 \}}
\frac{ a}{\tilde\l_{\rc,s-1}\tilde\e^\infty_{\rc,s}(1-\e_s)}\cdot$$
$$\int_{F_{i_1\dots i_{s-1}}\circ\g_j([0,1])}
||\tilde P_x\nabla u(x)||\cdot||\tilde P_x\nabla u(x)||\dr
(F_{i_1\dots i_{s-1}}\circ\g_j)_\sharp[
||DF_{i_1\dots i_{s-1}}\dot\g_j(t)||^2\dr t
]  =  $$
$$\int_G
||\tilde P_x\nabla u(x)||\cdot||\tilde P_x\nabla u(x)||\dr m(x)  .  \eqno (6.12)$$
An argument akin to that of lemma 6.1 implies that $m$ is a finite measure. 

The last formula, together with (6.2) and (6.10) implies that 
$$\ec_\rc(\nabla u,\nabla v)=
\int_G(T_x\nabla u(x),\nabla v(x))\dr\kappa(x)+$$
$$\int_G
(\tilde P_x\nabla u(x),\nabla v(x))\dr m(x)  .  \eqno (6.13)$$

Now note that the union of the vertices of all the pre-fractals $G_l$ is a subset of 
$\tilde  G$; being countable, it is null for Kusuoka's measure $\kappa$, which is non-atomic by [4]. On the other side, the set of vertices contains the ends of the cables, i. e. of  the segments $F_{i_1\dots i_{s-1}}\circ\g_j([0,1])$, with $j\in\Omega^s$; by (6.12), it is a null set also for $m$. In other words, the intersection of the support of $\kappa$ (which is the set $\tilde G$ of (1.24)) and that of the support of $m$ (which is the union of the "cables" 
$F_{i_1\dots i_{r}}\circ\g_j([0,1])$) in null for the measure $\mu$ defined by 
$$\mu=\kappa+m   .   \eqno (6.14)$$
Thus, if we define the following matrix field  
$$\tilde T_x=\left\{
\eqalign{
T_x&\txt{if}x\in\tilde G\cr
\tilde P_x&\txt{if}x\in G\setminus\tilde G          
}
\right. $$
then we can re-write formula (6.13) as 
$$\ec_\rc(\nabla u,\nabla v)=
\int_G(\tilde T_x\nabla u(x),\nabla v(x))\dr\mu(x)  .  \eqno (6.15)$$

We saw above that $T_x$ is bounded, while $\tilde P_x$ is bounded by definition; at the end, the Borel field of matrices $\tilde T_x$ is bounded. 

\vskip 2pc

\centerline{\bf \S 7}

\centerline{\bf The form is Dirichlet}

\vskip 1pc

\lem{7.1} Let $\fun{u}{\R^2}{\R}$ be affine and let $\ec_\rc$ be as in (6.10) or, more explicitly, as in (6.15). Then, $u$ is harmonic for $\ec$, i. e. formula (12) holds. 

\proof Let $v\in C^1(\R^2,\R)$ be such that $v(A)=v(B)=v(C)=0$; lemmas 3.2 and 5.1 imply that, for all $l\ge 1$, 
$$\ec_{\rc,l}(\nabla u,\nabla v)=0  .  $$
Now (12) follows from (6.11). 

\fin

\noindent{\bf Proof of theorem 1.} Point 1) is formula (6.11). 

Point 2) is lemma 7.1. 

We only sketch from [6] the proof of point 3) for the measure $\mu$ of (6.14): it is based on the standard calculation of the Laplacian of a composition. A side result will be Teplyaev's formula for the Laplacian, namely that, for all $\phi\in C^2(\R^2,\R)$  we have 
$$\Delta\phi(x)={\rm tr}(\tilde T_x D^2\phi(x))   
\txt{for $\mu$-a. e. $x\in G$. }   \eqno (7.1)$$

Let $\phi\in C^2(\R^2,\R)$ and let the test function $v\in C^1(\R^2,\R)$ satisfy $v(A)=v(B)=v(C)=0$. The first equality below is (6.15), the second one follows from a dumb application of the chain rule; the third equality needs no explanation and the fourth one comes from the Leibnitz formula for $C^1$ functions; the fifth one comes from the fact that the coordinate functions are affine and thus harmonic; as for the last equality, a simple computation shows that the integrands coincide: the expression in the round parentheses is the matrix whose two columns are 
$\nabla\partial_1\phi$ and 
$\nabla\partial_2\phi$. 
$$\ec_\rc(\nabla\phi,\nabla v)=
\int_G(\tilde T_x\cdot\nabla\phi,\nabla v)\dr\mu(x)=$$
$$\int_G(
\tilde T_x\cdot\partial_1\phi\cdot\nabla(x_1)+\tilde T_x\cdot\partial_2\phi\cdot\nabla(x_2),
\nabla v
)\dr\mu(x)=$$
$$\int_G(\tilde T_x\cdot\nabla(x_1),\partial_1\phi\cdot\nabla v)\dr\mu(x)+
\int_G(\tilde T_x\cdot\nabla(x_2),\partial_2\phi\cdot\nabla v)\dr\mu(x)=$$
$$\int_G(\tilde T_x\cdot\nabla(x_1),\nabla(\partial_1\phi\cdot v))\dr\mu(x)+
\int_G(\tilde T_x\cdot\nabla(x_2),\nabla(\partial_2\phi\cdot v))\dr\mu(x)-$$
$$\int_G(\tilde T_x\cdot\nabla(x_1),v\cdot\nabla\partial_1\phi)\dr\mu(x)-
\int_G(\tilde T_x\cdot\nabla(x_2),v\cdot\nabla\partial_2\phi)\dr\mu(x)=$$
$$-\int_G(\tilde T_x\cdot\nabla(x_1),\nabla\partial_1\phi)\cdot v\dr\mu(x)-
\int_G(\tilde T_x\cdot\nabla(x_2),\nabla\partial_2\phi)\cdot v\dr\mu(x)=$$
$$-\int_G{\rm tr}[
\tilde T_x\cdot(\nabla\partial_1(x),\nabla\partial_2\phi)
]\cdot v\dr\mu(x)  .  $$
Now (7.1) follows from the last formula. 

We show point 4), i. e. that $\ec_\rc$ extends to a local, regular Dirichlet form. It is local because it is the extension of local forms ([7]); regularity follows easily from the fact that 
$\ec$ is defined on $C^1(\R^2,\R)$ which is dense in $L^2(G,\mu)$ and $C(G,\R)$. We show closability. 

Actually, we prove the equivalent fact ([14]) that the quadratic form 
$\fun{}{u}{\ec(\nabla u,\nabla u)}$ is lower semicontinuous in $L^2(G,\mu)$: if 
$u_n,u\in C^1(\R^2,\R)$ and $u_n\tends u$ in $L^2(G,\mu)$, then 
$$\ec_\rc(\nabla u,\nabla u)\le\liminf_{n\tends+\infty}\ec_\rc(\nabla u_n,\nabla u_n)  .  $$
We sketch the proof of [6]. It follows easily from (6.15) that 
$$\ec_\rc(\nabla u,\nabla u)=\sup\ec_\rc(\nabla u,\phi\nabla v)$$
where the supremum is over all couples $\phi,v\in C^2(\R^2,\R)$ such that 

\noindent $\bullet$) $\phi(A)=\phi(B)=\phi(C)=0$ and 

\noindent $\bullet$) 
$\ec_\rc(\phi\nabla v,\phi\nabla v)\le \ec_\rc(\nabla u,\nabla u)$. 

Thus, the thesis follows if we prove that the function from $C^1(\R^2,\R)$ to $\R$ 
$$\fun{}{u}{\ec_\rc(\nabla u,\phi\nabla v)}$$
is continuous for the $L^2(G,\mu)$ topology for all couples $(\phi,v)$ as above. We show this: the first equality below follows from (6.15); the second one from the Leibnitz formula for $C^1$ functions and (6.15) again; the third one follows from the integration by parts formula. 
$$\ec_\rc(\nabla u,\phi\nabla v)=
\ec_\rc(\phi\nabla u,\nabla v)=$$
$$\ec_\rc(\nabla(\phi u),\nabla v)-\ec_\rc(u\nabla\phi,\nabla v)=$$
$$-\int_G\Delta v\cdot\phi\cdot u\dr\mu-\ec_\rc(u\nabla\phi,\nabla v)  .  $$
Since $v\in C^2(\R^2,\R)$ and $\tilde T_x$ is bounded, (7.1) implies that $\Delta v$ is bounded too; since $\phi$ is bounded and $\mu$ is a finite measure, we get that the function 
$$\fun{}{u}{
\int_G\Delta v\cdot\phi\cdot u\dr\mu
}$$
is continuous from $L^2(G,\mu)$ to $\R$. The fact that 
$$\fun{}{u}{
\ec(u\nabla\phi,\nabla v)
}$$
is continuous follows from (6.15) and the fact that $\tilde T_x$, $\nabla\phi$ and $\nabla v$ are bounded on the compact set $G$. 

As for point 5), it follows taking limits in the recurrence (3.6) and recalling (5.1) and (6.11).

\fin

\vskip 2pc
\centerline{\bf References} 






\noindent [1] P. Alonso-Ruiz, U. Freiberg, J. Kigami, Completely symmetric resistance forms on the Stretched Sierpinski Gasket, Journal of Fractal Geometry, {\bf 5-3}, 227-277, 2018. 

\noindent [2] P. Alonso Ruiz, D. J. Kelleher, A. Teplyaev, Energy and Laplacian on Hanoi-type quantum graphs, Journal of Physics A: Mathematical and Theoretical, {\bf 49}, 165-206, 2016. 

\noindent [3] A. Arauza Rivera, A. Lankford, M. McClinton, S. Torres, An IFS for the stretched level-$n$ Sierpinski gasket, PUMP Journal of Undergraduate Research, {\bf 5}, 105-121, 2022.  

\noindent [4] U.Bessi, Another point of view on Kusuoka's measure, DCDS, {\bf 41-7}, 3251-3271, 2021.  

\noindent [5] U. Bessi, Families  of Kusuoka-like measures on fractals, preprint. 

\noindent [6] U. Bessi, Total variation of smooth functions on harmonic fractals, preprint. 












\noindent [7] M. Fukushima, Y. Oshima, M. Takeda, Dirichlet forms and symmetric Markov processes, De Gruyter, G\"ottingen, 2011. 


\noindent [8] A. Johansson, A. \"Oberg, M. Pollicott, Ergodic theory of Kusuoka's measures, J. Fractal Geom., {\bf 4}, 185-214, 2017.  

\noindent [9] N. Kajino, Analysis and geometry of the measurable Riemannian structure on the Sierpiski gasket, Contemporary Mathematics, {\bf 600}, 91-133, 2013. 


\noindent [10] J. Kigami, Analysis on fractals, Cambridge, 2001. 

\noindent [11] J.Kigami, Harmonic Analysis for Resistance Forms, Journal of Functional Analysis, {\bf 204-2}, 399-444, 2003. 



\noindent [12] S. Kusuoka, Dirichlet forms on fractals and products of random matrices, Publ. Res. Inst. Math. Sci., {\bf 25}, 659-680, 1989. 





\noindent [13] I. D. Morris, Ergodic properties of matrix equilibrium state, Ergodic Theory and Dyn. Sys., {38/6}, 2295-2320, 2018.

\noindent [14] U.Mosco, Composite media and asymptotic Dirichlet forms, J. Functional Analysis, {\bf 123}, 368-421, 1994.

\noindent [15] R. Peirone, Convergence of Dirichlet forms on fractals, mimeographed notes. 



\noindent [16] W. Perry, M. Pollicott, Zeta functions and the periodic orbit structure of hyperbolic dynamics, Asterisque, {\bf 187-188}, 1990. 

\noindent [17] M. Piraino, The weak Bernoulli property for matrix equilibrium states, preprint on Arxiv.


\noindent [18] W. Rudin, Real and complex Analysis, New Delhi, 1983. 




\noindent [19] A. Teplyaev, Energy and Laplacian on the Sierpinski gasket, Proceedings of Symposia in Pure Mathematics, {\bf 72-1}, 131-154, Providence, R. I., 2004. 


\noindent [20] M. Viana, Stochastic analysis of deterministic systems, mimeographed notes.



\end